\documentclass[journal,twoside,web]{ieeetran}
\usepackage{amsthm}

\usepackage{cite}
\usepackage{amsmath,amssymb,amsfonts}
\usepackage{algorithm}
\usepackage{algpseudocode}
\usepackage{bm}
\usepackage{breqn}
\newtheorem{theorem}{Theorem}

\newtheorem{lemma}{Lemma}
\newtheorem{definition}{Definition}
\newtheorem{assumption}{Assumption}
\newtheorem{proposition}{Proposition}
\newtheorem{cond}{Condition}

\usepackage{caption}
\usepackage{enumerate}
\usepackage[caption=false, font=footnotesize]{subfig}

\allowdisplaybreaks

\begin{document}

\title{\huge Decentralized State-Dependent Markov Chain Synthesis \\ with an Application to Swarm Guidance}

\author{Samet Uzun, \IEEEmembership{Student Member, IEEE}, Nazım Kemal Üre, \IEEEmembership{Member, IEEE}, Behçet Açıkmeşe, \IEEEmembership{Fellow, IEEE}
\thanks{ Samet Uzun and Behçet Açıkmeşe are with the William E. Boeing Department of Aeronautics and Astronautics, University of Washington, Seattle, WA 98125, USA (e-mail: samet@uw.edu; behcet@uw.edu). }
\thanks{Nazım Kemal Üre is with the 
Department of Artificial Intelligence and Data Engineering,
Istanbul Technical University, Istanbul, 34469, Turkey (e-mail: ure@itu.edu.tr).}
\thanks{\textcopyright 2024 IEEE. Personal use of this material is permitted. Permission from IEEE must be obtained for all other uses, in any current or future media, including reprinting/republishing this material for advertising or promotional purposes, creating new collective works, for resale or redistribution to servers or lists, or reuse of any copyrighted component of this work in other works.}}

\maketitle
\begin{abstract}
This paper introduces a decentralized state-dependent Markov chain synthesis (DSMC) algorithm for finite-state Markov chains. 
We present a state-dependent consensus protocol that achieves exponential convergence under mild technical conditions, without relying on any connectivity assumptions regarding the dynamic network topology.
Utilizing the proposed consensus protocol, we develop the DSMC algorithm, updating the Markov matrix based on the current state while ensuring the convergence conditions of the consensus protocol.
This result establishes the desired steady-state distribution for the resulting Markov chain, ensuring exponential convergence from all initial distributions while adhering to transition constraints and minimizing state transitions.
The DSMC's performance is demonstrated through a probabilistic swarm guidance example, which interprets the spatial distribution of a swarm comprising a large number of mobile agents as a probability distribution and
utilizes the Markov chain to compute transition probabilities between states.
Simulation results demonstrate faster convergence for the DSMC based algorithm when compared to the previous Markov chain based swarm guidance algorithms.
\end{abstract}

\begin{IEEEkeywords}
Markov Chains, Consensus Protocol, Decentralized Control, 
Probabilistic Swarm Guidance.
\end{IEEEkeywords}

\section{Introduction} \label{sec:intro}

Markov chain synthesis has garnered attention from various disciplines, including physics, systems theory, computer science, and numerous other fields of science and engineering. This attention is particularly notable within the context of Monte Carlo Markov Chain (MCMC) algorithms \cite{levin2017markov, fishman2013monte, gilks1995markov}.
The fundamental idea underlying MCMC algorithms is to synthesize a Markov chain that converges to a specified steady-state distribution. 
Random sampling of a large state space while adhering to a predefined probability distribution is the predominant use of MCMC algorithms.
The current literature covers a broad spectrum of methodologies for Markov chain synthesis, incorporating both heuristic approaches and optimization-based techniques \cite{hastings1970monte, boyd2004fastest, boyd2009fastest}. Each method provides specialized algorithms tailored to the synthesis of Markov chains in alignment with specific objectives or constraints.
Markov chain synthesis plays a central role in probabilistic swarm guidance, which has led to the development of various algorithms incorporating additional transition and safety constraints \cite{accikmecse2012markov, el2018safe, accikmecse2015markov, demir2014density, accikmecse2015convex, demir2015decentralized, demirer2017safe, bandyopadhyay2017probabilistic, jang2018local, djeumou2020probabilistic, djeumou2022probabilistic}. 
The probabilistic swarm guidance algorithm models the spatial distribution of the swarm
agents as a probability distribution and employs the synthesized Markov chain to guide the spatial distribution of the swarm.

\vspace{0.1cm}
Consensus protocols form an important field of research that has a strong connection with Markov chains \cite{mesbahi2010graph}. 
Consensus protocols are a set of rules used in distributed systems to achieve agreement among a group of agents on the value of a variable \cite{olfati2007consensus, cao2012overview, qin2016recent, kia2019tutorial}.
Markov chains are often employed to model and analyze the dynamics and convergence properties of consensus protocols, providing insights into their behavior and performance \cite{xiao2004fast, boyd2006randomized, tahbaz2007small, olshevsky2013degree}.

\vspace{0.1cm}
The current paper presents a consensus protocol specifically tailored to operate on a dynamic graph topology. We establish that the protocol provides exponential convergence guarantees, even under mild technical conditions. Consequently, our approach eliminates the reliance on conventional connectivity assumptions commonly found in the existing literature \cite{jadbabaie2003coordination, moreau2005stability, ren2005consensus, fax2004information,  olfati2004consensus, sluvciak2016consensus}.
Building on this new consensus protocol, the paper introduces a decentralized state-dependent Markov chain (DSMC) synthesis algorithm. It is demonstrated that the synthesized Markov chain, formulated using the proposed consensus algorithm, satisfies the aforementioned mild conditions. This, in turn, ensures the exponential convergence of the Markov chain to the desired steady-state.
Subsequently, the DSMC algorithm's performance is demonstrated on the probabilistic swarm guidance problem, with the objective of controlling the spatial distribution of swarm agents. Through simulations, it is shown that the DSMC algorithm achieves considerably faster convergence
compared to other existing Markov chain synthesis algorithms.

\subsection{Related Works} \label{sec:related_work}

The Metropolis-Hastings algorithm is widely recognized as one of the most prominent techniques for MCMC algorithms \cite{hastings1970monte}.
For the fastest mixing Markov chain synthesis, the problem is formulated as a convex optimization problem in \cite{boyd2004fastest}, assuming that the Markov chain is symmetric. This paper also presents an extension to the method that involves synthesizing the fastest mixing reversible Markov chain with a given desired distribution. Furthermore, the number of variables in the optimization problem is reduced in \cite{boyd2009fastest} by exploiting the symmetries in the graph. 

\vspace{0.1cm}
The probabilistic swarm guidance problem has been a crucial domain for the application and enhancement of Markov chain synthesis algorithms.
A comprehensive review of the broader category of multi-agent algorithms is presented in \cite{rossi2018review}, while a survey specifically focusing on aerial swarm robotics is provided in \cite{chung2018survey}. Additionally, \cite{schranz2020swarm} offers an overview of existing swarm robotic applications.
For swarm guidance purposes, certain deterministic algorithms have been developed in \cite{lumelsky1997decentralized, richards2002spacecraft, tillerson2002co, scharf2003survey, kim2004multiple, ramirez2010distributed}. However, these algorithms may become computationally infeasible when dealing with swarms that comprise hundreds to thousands of agents.
In the context of addressing the guidance problem for a large number of agents, considering the spatial distribution of swarm agents and directing it towards a desired steady-state distribution offers a computationally efficient approach. In this regard, both probabilistic and deterministic swarm guidance algorithms are presented in \cite{zhao2011density, eren2017velocity, elamvazhuthi2020controllability, elamvazhuthi2019mean, biswal2021decentralized, krishnan2018distributed1, krishnan2018distributed2} for continuous state spaces.
For discrete state spaces, a probabilistic guidance algorithm is introduced in \cite{accikmecse2012markov}. 
This algorithm treats the spatial distribution of swarm agents, called the density distribution, as a probability distribution and employs the Metropolis-Hastings (M-H) algorithm to synthesize a Markov chain that guides the density distribution toward a desired state.
The probabilistic guidance algorithm led to the development of numerous Markov chain synthesis algorithms involving specific objectives and constraints \cite{el2018safe, accikmecse2015markov, demir2014density, accikmecse2015convex, demir2015decentralized, demirer2017safe, bandyopadhyay2017probabilistic, jang2018local, djeumou2020probabilistic, djeumou2022probabilistic}.
In \cite{el2018safe}, the Metropolis-Hastings algorithm is extended to incorporate safety upper bound constraints on the probability vector. This paper includes numerical simulations that demonstrate the application of the extension in a probabilistic swarm guidance problem.
In order to enhance convergence rates, \cite{accikmecse2015markov} introduces a convex optimization-based technique for Markov chain synthesis. This technique formulates the objective function and constraints using linear matrix inequalities (LMIs) to synthesize a Markov chain capable of achieving a desired distribution while adhering to specified transition constraints. Notably, this study does not impose any assumptions on the Markov chain, rendering the problem inherently non-convex. The problem is convexified for practical purposes but the optimal convergence rate of the original non-convex problem cannot be attained.
In \cite{demir2014density}, the approach presented in \cite{accikmecse2015markov} is enhanced by incorporating state-feedback to further improve the convergence rate. 
These works are also extended to impose density upper bounds and density rate constraints in \cite{accikmecse2015convex} and density flow and density diffusivity constraints in \cite{demir2015decentralized}.
In \cite{demirer2017safe}, a more general approach is proposed for cases where the environment is stochastic, and the resulting Markov chain is shown to be a linear function of the stochastic environment and decision policy.
These convex optimization problems can be solved using the interior-point algorithm but the computation time of the algorithm increases rapidly with increasing dimensionality of the state space \cite{boyd1994linear}.
Furthermore,
neither the M-H algorithm nor convex optimization-based approaches 
attempt to minimize the number of state transitions.
The swarm guidance problem is formulated as an optimal transport problem in \cite{bandyopadhyay2014probabilistic} to minimize the number of agent transitions. However, besides the similar computational complexity issues,
the performance of the algorithm drops significantly if the current density distribution of the swarm cannot be estimated accurately.
The time-inhomogeneous Markov chain approach to the probabilistic swarm guidance problem (PSG-IMC algorithm) is developed in \cite{bandyopadhyay2017probabilistic} to minimize the number of state transitions. This algorithm is computationally efficient and yields reasonable results with low estimation errors. 
However, all feedback-based algorithms mentioned above require global feedback on the state of the density distribution.
Communication between all agents has to be established to estimate the density distribution in the probabilistic swarm guidance problem. 
Generating perfect feedback of the density distribution is often impractical, leading to the routine occurrence of estimation errors.
An alternative approach that requires only local information is developed in \cite{jang2018local}. In this work, the time-inhomogeneous Markov chain approach presented in \cite{bandyopadhyay2017probabilistic} is modified to work with local information, and the algorithm is used for minimizing the number of state transitions. 
Nevertheless, in both global and local information based time-inhomogeneous Markov chain approaches, it is observed that the convergence rate diminishes significantly when the state transition capabilities become more restricted.
As it can be seen in Corollary 1 in \cite{bandyopadhyay2017probabilistic} or Corollary 3 in \cite{jang2018local}, the transition probability from a state to any directly connected state cannot be higher than the desired density value of the corresponding directly connected state.
In situations where there are sparsely connected regions in the state space, it is common to observe a relatively low sum of desired density values among directly connected states.
Consequently, there is a higher probability of remaining in the same state rather than transitioning to other states.
In terms of the convergence rate, these algorithms are only effective in cases with high transition capabilities.
Additionally, the performance of these algorithms is highly sensitive to hyperparameters and requires careful selection for optimum results in each experiment.

\vspace{0.1cm}
Graph temporal logic (GTL) is introduced in \cite{djeumou2020probabilistic} to impose high-level task specifications as a constraint to the Markov chain synthesis. Markov chain synthesis is formulated as mixed-integer nonlinear programming (MINLP) feasibility problem and the problem is solved using a coordinate descent algorithm. In addition, an equivalence is proven between the feasibility of the MINLP and the feasibility of a mixed-integer linear program (MILP) for a particular case where the agents move along the nodes of a complete graph. While this study assumes homogeneous swarms for Markov chain synthesis subject to finite-horizon GTL formulas, an improved version of the formulation is presented in \cite{djeumou2022probabilistic} to enable probabilistic control of heterogeneous swarms subject to infinite-horizon GTL formulas. Instead of solving the resulting MINLP using a coordinate descent algorithm, a sequential scheme, which is faster, more accurate, and robust to the choice of the starting point, is developed in the aforementioned paper.

\vspace{0.1cm}
Markov chains and consensus protocols share a close relationship. The rich theory of Markov chains has proven to be valuable in analyzing specific consensus protocols. Notable works such as \cite{xiao2004fast, boyd2006randomized, tahbaz2007small, olshevsky2013degree} have leveraged Markov chain theory to provide insights and analysis for consensus protocols.
Consensus protocols, in contrast to Markov chains, operate without the limitations of non-negative nodes and edges or the requirement for the sum of nodes to equal one \cite{mesbahi2010graph}. This broader scope enables consensus protocols to address a significantly wider range of problem spaces.
Therefore, there is a significant interest in consensus protocols in a broad range of multi-agent networked systems research, including distributed coordination of mobile autonomous agents \cite{jadbabaie2003coordination, moreau2005stability, ren2005consensus, fax2004information,  olfati2004consensus, oh2015survey}, distributed optimization \cite{nedic2009distributed, boyd2011distributed, shi2015extra, kia2015distributed, yang2019survey}, distributed state estimation \cite{olfati2007distributed, kamal2013information}, or dynamic load-balancing for parallel processors \cite{cybenko1989dynamic, xiao2007distributed}.
There are comprehensive survey papers that review the research on consensus protocols  \cite{olfati2007consensus, cao2012overview, qin2016recent, kia2019tutorial}. In many scenarios, the network topology of the consensus protocol is a switching topology due to failures, formation reconfiguration, or state-dependence. There is a large number of papers that propose consensus protocols with switching network topologies and convergence proofs of these algorithms are provided under various assumptions \cite{jadbabaie2003coordination, moreau2005stability, ren2005consensus, fax2004information,  olfati2004consensus, sluvciak2016consensus}.
In \cite{jadbabaie2003coordination}, a consensus protocol is proposed to solve the alignment problem of mobile agents, where the switching topology is assumed to be periodically connected.
This assumption means that
the union of networks over a finite time interval is strongly connected. 
Another algorithm is proposed in \cite{moreau2005stability} that assumes the underlying switching network topology is ultimately connected. This assumption means that the union of graphs over an infinite interval is strongly connected. In \cite{ren2005consensus}, previous works are extended to solve the consensus problem on networks under limited and unreliable information exchange with dynamically changing interaction topologies. The convergence of the algorithm is provided under the ultimately connected assumption. 
Another consensus protocol is introduced in \cite{fax2004information} for the cooperation of vehicles performing a shared task using inter-vehicle communication. Based on this work, a theoretical framework is presented in \cite{olfati2004consensus} to solve consensus problems under a variety of assumptions on the network topology such as strongly connected switching topology.
In \cite{sluvciak2016consensus}, a consensus protocol with state-dependent weights is proposed and it is 
assumed that corresponding graphs are weakly connected, which means that
the network is assumed to be connected over the iterations. 
Additionally, optimization-based algorithms are proposed to obtain a high convergence rate for the consensus protocol in \cite{xiao2004fast} and \cite{kim2005maximizing}.

\subsection{Main Contributions}
In this paper, 
we first propose a consensus protocol with state-dependent weights.
The proposed consensus protocol does not require any connectivity assumption on the dynamic network
topology, unlike the existing methods in the literature \cite{jadbabaie2003coordination, moreau2005stability, ren2005consensus, fax2004information,  olfati2004consensus, sluvciak2016consensus}.
We provide an exponential convergence guarantee for the consensus protocol under some mild technical conditions, which can be verified straightforwardly.
We then present a decentralized Markov-chain synthesis (DSMC) algorithm based on the proposed consensus protocol and we prove that the resulting DSMC algorithm satisfies these mild conditions.
This result is employed to prove that the resulting Markov chain has a desired steady-state distribution and that all initial distributions exponentially converge to this steady-state.
Unlike the homogeneous Markov chain synthesis algorithms in \cite{hastings1970monte, accikmecse2012markov, boyd2004fastest,
boyd2009fastest, el2018safe, accikmecse2015markov}, the Markov matrix, synthesized by our algorithm, approaches the identity matrix as the probability distribution converges to the desired steady-state distribution. Hence the proposed algorithm attempts to minimize the number of state transitions, which eventually converge to zero as the probability distribution converges to the desired steady-state distribution.
Whereas previous time-inhomogeneous Markov chain synthesis algorithms in \cite{bandyopadhyay2017probabilistic, jang2018local} only provide asymptotic convergence, the DSMC algorithm provides an exponential convergence rate guarantee.
Furthermore, unlike previous algorithms in \cite{bandyopadhyay2017probabilistic, jang2018local}, the convergence rate of the DSMC algorithm does not rapidly decrease in scenarios where the state space contains sparsely connected regions.
Due to the decentralized nature of the consensus protocol, the Markov chain synthesis relies on local information, similar to the approach presented in \cite{jang2018local},
and a complex communication architecture is not required for the estimation of the distribution. 
By presenting numerical evidence within the context of the probabilistic swarm guidance problem, we demonstrate that the convergence rate of the swarm distribution to the desired steady-state distribution is substantially faster when compared to previous methodologies.
In summary: 
\begin{itemize}
    \item We propose a consensus protocol with state-dependent weights and prove its exponential convergence under mild technical conditions, without relying on the typical connectivity assumptions associated with the underlying graph topology.
    \item Based on the proposed consensus protocol, we introduce the DSMC algorithm, which is shown to meet the convergence conditions outlined by the consensus protocol.
    \item  
    Simulation results demonstrate that the DSMC algorithm achieves faster convergence, characterized by an exponential convergence guarantee, compared to existing homogeneous and time-inhomogeneous Markov chain synthesis algorithms presented in \cite{accikmecse2012markov} and \cite{bandyopadhyay2017probabilistic}.
\end{itemize}

\vspace{-0.2cm}
\subsection{Organization and Notation}
The paper is organized as follows. Section \ref{sec:DSDN} presents the consensus protocol with state-dependent weights. The decentralized state-dependent Markov matrix synthesis (DSMC) algorithm is introduced in Section \ref{sec:somm}. 
Section \ref{sec:background} introduces the probabilistic swarm guidance problem formulation, and presents numerical simulations of swarms converging to desired distributions. The paper is concluded in Section \ref{sec:conclusion}.

\textit{Notation:}
$\mathbb{R}$ and $\mathbb{C}$ represent the set of real numbers and complex numbers, respectively. $\mathbb{R}_+$ and  $\mathbb{Z}_+$ represent the set of non-negative real numbers and non-negative integers, respectively. $\mathbb{R}^n$ is the $n$ dimensional real vector space.
$\bm{0}$ is a zero matrix and $\bm{1}$ is a matrix of ones in appropriate dimensions. 
$x[i](k)$ represents the $i^{\text{th}}$ element of the vector $x \in \mathbb{R}^n$ at time $k$.
$\left\| x \right\|_p$ denotes the $\ell_p$ vector norm. 
$(v_1, v_2, ... , v_n)$ represents a vector, such that $(v_1, v_2, ... , v_n) \equiv [v_1^T, v_2^T, ... , v_n^T]^T$ where $v_i$ have arbitrary dimensions. 
$M > (\geq) \; \bm{0}$ implies that $M \in \mathbb{R}^{m \times n}$ is a positive (non-negative) matrix. 
$\mathcal{I}_1^m$ denotes the set of integers such that $\mathcal{I}_1^m = \{ 1, 2, \dots,  m \}$.
$\lambda_i(A)$ is the $i^{th}$ eigenvalue of a matrix $A \in \mathbb{R}^{m \times m}$ such that $\lambda_i(A) \leq \lambda_{i+1}(A)$ for $i \in \mathcal{I}_1^{m-1}$.
$\sigma (A)$ is the set of eigenvalues of $A \in \mathbb{R}^{m \times m}$ such that $\sigma (A) = \{ \lambda_1, \lambda_2, \dots,  \lambda_m \}$.
$\emptyset$ denotes the empty set.
$V \setminus W$ represents the elements in set $V$ that are not in set $W$. 
$\mathcal{P}$ denotes the probability of a random variable. 

\newpage
\section{A Consensus Protocol \\ with State-Dependent Weights} \label{sec:DSDN}
\vspace{-0.12cm}
The consensus protocol is a process used for achieving an agreement among distributed agents. 
In this section,
we introduce a consensus protocol with state-dependent weights to reach a consensus on time-varying weighted graphs.
Unlike other proposed consensus protocols in the literature, the consensus protocol we introduce does not require any connectivity assumption on the dynamic network topology. We provide theoretical analysis for proof of exponential convergence under some mild technical conditions.
We then use the proposed consensus protocol for the synthesis of the column stochastic Markov matrix in Section \ref{sec:somm}. It is proven that these mild assumptions required by the proposed consensus protocol are satisfied by the Markov matrix synthesis algorithm.
We first define the graph for our consensus approach.
\vspace{-0.2cm}
\begin{definition} \label{def:graph_def}
    (\textit{Time-varying weighted graph}) A time-varying weighted graph is a tuple, $\mathcal{G}_w(k) = (\mathcal{V}, \mathcal{E}, w(k))$, where  $\mathcal{V} = \{ \mathrm{v}_1, ..., \mathrm{v}_m \}$ is the set of vertices, $\mathcal{E} \subseteq \mathcal{V} \times \mathcal{V}$ is the set of undirected edges, and $w(k)$ is a time-varying function such that $w : \mathcal{E} \xrightarrow[]{} \mathbb{R}_{+}$,
    where $w_{i,j}(k)$ represents the value of weight for edge $ \{ \mathrm{v}_i,  \mathrm{v}_j \}  \in \mathcal{E}$ at time $k$. 
    \begin{itemize}
        \item (\textit{Values of vertices}) 
         A time-varying vector $e(k) \in \mathbb{R}^m$ such that $ \bm{1}^T e(k) = 0$ for $k \in \mathbb{Z}_+$ can define a vector of values of vertices, i.e.,   $e[i](k)$ represents the value of the vertex $\mathrm{v}_i$ at time $k$.
        
        \item (\textit{Uniform graph induced by time-varying weighted graph}) 
       $\mathcal{G} = (\mathcal{V}, \mathcal{E})$ is the uniform graph induced by $\mathcal{G}_w(k)$  and it is obtained by setting $w_{i,j}(k) = 1 $, where $\{\mathrm{v}_i, \mathrm{v}_j\} \in \mathcal{E}$ for all $k \in \mathbb{Z}_+$.
 
        \item (\textit{Adjacency matrices}) 
        Two vertices $\mathrm{v}_i$ and $\mathrm{v}_j$ are called adjacent if $ \{ \mathrm{v}_i,  \mathrm{v}_j \}  \in  \mathcal{E}$. 
        The adjacency matrix of the uniform graph $\mathcal{G}$ is represented by $A_a(\mathcal{G})$, where $A_a(\mathcal{G})[i,j] = 1$, if $\mathrm{v}_i$ and $\mathrm{v}_j$ are adjacent, and is $0$ otherwise. The adjacency matrix of the time-varying weighted graph $\mathcal{G}_w(k)$ is represented by $A_a(\mathcal{G}_w(k))$, where $A_{a}(\mathcal{G}_w(k))[i,j] = A_a(\mathcal{G})[i,j] w_{i,j}(k)$.
        
        \item (\textit{Degree matrices}) 
        $\mathcal{D}(\mathcal{G})$ and $\mathcal{D}(\mathcal{G}_w(k))$ are the diagonal degree matrices of the graphs $\mathcal{G}$ and $\mathcal{G}_w(k)$ such that $\mathcal{D}(\mathcal{G})[i,i] = \sum_j A_a(\mathcal{G})[i,j]$ and $\mathcal{D}(\mathcal{G}_w(k))[i,i] = \sum_{j \in \mathcal{I}_1^m } A_{a}(\mathcal{G}_w(k))[i,j]$.
        
        \item (\textit{Laplacian matrices}) 
        $\mathcal{L}(\mathcal{G})$ and $\mathcal{L}(\mathcal{G}_w(k))$ are the Laplacian matrices of the graphs $\mathcal{G}$ and $\mathcal{G}_w(k)$ such that $\mathcal{L}(\mathcal{G}) = \mathcal{D}(\mathcal{G}) - A_a(\mathcal{G})$ and $\mathcal{L}(\mathcal{G}_w(k)) = \mathcal{D}(\mathcal{G}_w(k)) - A_{a}(\mathcal{G}_w(k))$.
        
        \item  (\textit{Activated and Deactivated Edges}) If $w_{i,j}(k) = 0$, where $\{ \mathrm{v}_i, \mathrm{v}_j\} \in \mathcal{E}$, the edge $\{ \mathrm{v}_i, \mathrm{v}_j\}$ is  deactivated. $\mathcal{E}_A(k)$ and $\mathcal{E}_D(k)$ represent the activated and deactivated edges for time  $k$ such that $\mathcal{E} = \mathcal{E}_A(k) \cup \mathcal{E}_D(k)$.
    \end{itemize}
\end{definition}
\vspace{-0.2cm} 
For convenience, we use $A_{a}$, $A_{a_w}(k)$, $\mathcal{D}$, $\mathcal{D}_w(k)$, $\mathcal{L}$ and $\mathcal{L}_w(k)$ instead of using $A_a(\mathcal{G})$, $A_a(\mathcal{G}_w(k))$, $\mathcal{D}(\mathcal{G})$, $\mathcal{D}(\mathcal{G}_w(k))$, $\mathcal{L}(\mathcal{G})$ and $\mathcal{L}(\mathcal{G}_w(k))$, respectively.

The following assumption on the topology of the uniform graph $\mathcal{G}$ is needed for convergence analysis.
\vspace{-0.2cm}
\begin{assumption} \label{asm:str_con_graph}
    The uniform graph $\mathcal{G}$ is a connected graph, which means there exists a path between all vertices, or equivalently, $(I + A_{a})^{m-1} > \bm{0}$ for $A_{a} = A_{a}^T \in \mathbb{R}^{m \times m}$ \cite[section 2.1]{mesbahi2010graph}.
\end{assumption}

We will consider time-varying weighted graphs where the weights of the edges depend on the values of vertices. The following definition is needed to present the relation between the weights of the edges and the values of vertices.
\begin{definition} \label{def:erridxset}
    \textit{(Index sets with respect to values of vertices)}
    The index sets $I_p(k)$ and $I_n(k)$ contain the indices of the non-negative and negative valued vertices for time $k$, i.e., 
    \vspace{-0.12cm}
    \begin{equation*}
        \begin{aligned}
            I_p(k) &= \{ i : e[i](k) \geq 0 \}, \\
            I_n(k) &= \{ i : e[i](k) < 0 \}. \\
        \end{aligned}
    \end{equation*}
    The index set $I_{n_p}(k) \subseteq I_{n}(k)$ consists of the indices of the negative valued vertices that are adjacent to the non-negative valued vertices, i.e., 
    \vspace{-0.12cm}
    \begin{equation*}
        \begin{aligned}
            I_{n_p}(k) &= \{ i : i \in I_n(k), j \in I_p(k) \text{ and } A_{a}[i,j] = 1 \}.
        \end{aligned}
    \end{equation*}
    The set that contains the edges between $\mathrm{v}_i$, $i \in I_n(k)$ and $\mathrm{v}_j$, $j \in I_{n_p}(k)$ is defined as
    \begin{equation*}
        \begin{aligned}
            I_{{n_p},n}(k) &= \{ \{ i, j \} : i \in I_n(k), j \in I_{n_p}(k) \text{ and } A_{a}[i,j] = 1\}.
        \end{aligned}
    \end{equation*}
\end{definition}

We provide a condition that represents a relation between the values of the vertices and the weights of the edges of the graph $\mathcal{G}_w(k)$. 
This condition is the key to providing the convergence proof of the consensus protocol with state-dependent weights.
\begin{cond} \label{con:pos_w}
    \hfill 
    \begin{enumerate}[\qquad(a)]
    \item If $i \in I_{p}(k)$, then there exists a $c_1$ such that $0 < c_1 \leq w_{i,j}(k)$ for all $\{ \mathrm{v}_i, \mathrm{v}_j \} \in \mathcal{E}$. \label{item:c1a}
    \item If $\{ i, j \} \in I_{{n_p}, n}(k)$, then $w_{i,j}(k) = 0$. \label{item:c1b}
    \end{enumerate}
\end{cond}
\vspace{0.15cm}

In other words, Condition \ref{con:pos_w}\ref{item:c1a} implies that the weight of an edge that connects a non-negative valued vertex to any other vertex cannot be arbitrarily close to $0$. Therefore, if $A_{a_w}[i,j](k) = 0$ where $e[i](k) \geq 0$, then $A_{a}[i,j](k) = 0$. Condition \ref{con:pos_w}\ref{item:c1b} implies that no transition is allowed between two negative valued vertices if at least one of them is adjacent to a non-negative valued vertex. 
Edges that have zero weights are called deactivated edges.
Note that although $\mathcal{G}$ is assumed to be a connected graph in Assumption \ref{asm:str_con_graph}, $\mathcal{G}_w(k)$ may not be a connected graph for some $k$, due to these deactivated edges.

\begin{figure}[!hbt]
    \centering
    \includegraphics[width=0.7\linewidth]{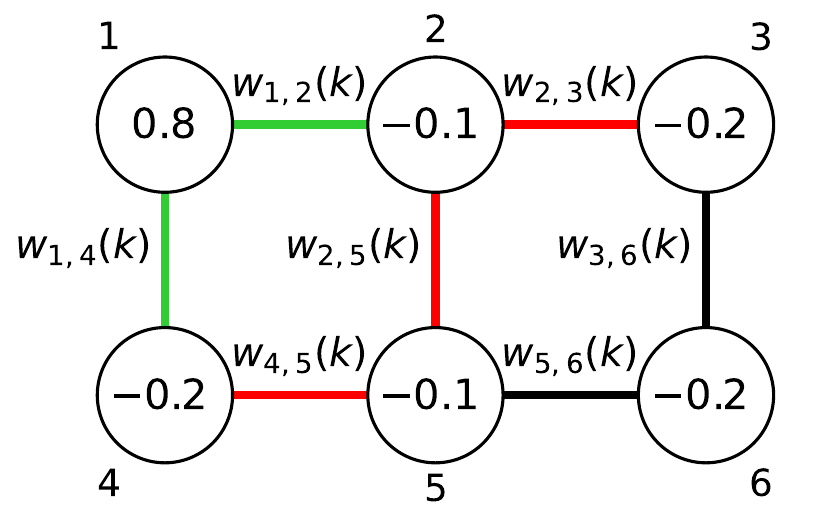}
    \caption{Vertices of the graph and corresponding values are represented by circles. The indices of the vertices are shown in their corners. Edges of the vertices are represented by lines and weights of the edges are represented on the sides of the edges. Green and red lines represent the edges that are implied by Condition \ref{con:pos_w}\ref{item:c1a} and Condition \ref{con:pos_w}\ref{item:c1b}, respectively. Black lines represent the remaining edges of the graph.}
    \label{fig:cond_graph_ex_g}
\end{figure}

The following example is provided to help the reader interpret the implication of Condition \ref{con:pos_w}.
\paragraph*{Example 1}
The values of vertices and weights of the edges are represented on a graph for a time $k$ in Figure \ref{fig:cond_graph_ex_g}. According to Definition \ref{def:erridxset}, $I_p(k) = \{ 1 \}$, $I_n(k) = \{ 2, 3, 4, 5, 6 \}$, $I_{n_p}(k) = \{ 2, 4 \}$ and $I_{{n_p},n}(k) = \{ \{ 2, 3 \} \{ 2, 5 \} \{ 4, 5 \} \}$. Condition \ref{con:pos_w}\ref{item:c1a} implies that $0 < c_1 \leq w_{1,2}(k)$ and $0 < c_1 \leq w_{1,4}(k)$ since $0 \leq e[1](k)$. Condition \ref{con:pos_w}\ref{item:c1b} implies that weights of all other edges of the nodes $\mathrm{v}_2$ and $\mathrm{v}_4$ are zero, which means that $w_{2,3}(k) = w_{2,5}(k) = w_{4,5}(k) = 0$.

The following lemma shows that each non-negative valued vertex is connected to a negative valued vertex in the graph $\mathcal{G}_w$ under Condition \ref{con:pos_w}.
\begin{lemma} \label{lem:err_vec}
    Assume that Condition \ref{con:pos_w} holds. If $e(k) \neq \bm{0}$, then for each $i$ such that $e[i](k) \geq 0$, $\exists j$ such that $e[j](k) < 0$, where $A_{a_w}^n[i,j](k)\neq 0$ for some exponent $n \in \mathbb{Z}_+$.
\end{lemma}
\begin{proof}
    Since $\bm{1}^T e(k) = 0$ by Definition \ref{def:graph_def} 
    and $e(k) \neq \bm{0}$, $\exists i,j$ such that $e[i](k) > 0$ and $e[j](k) < 0$. 
    
    The lemma is proven by contradiction.
    Suppose that there exists 
    one non-negative valued vertex $\mathrm{v}_i$ such that $e[i](k) \geq 0$ and for all $j $ such that $e[j](k) < 0$, $A_{a_w}^n[i,j](k)= 0$, $\forall n \in \mathbb{Z}_+$. In other words, there exists one non-negative valued vertex that is not connected to any negative valued vertices. 
    Denote the set of such non-negative valued vertices as $\mathcal{V}^+(k)$. The set of remaining vertices
    is denoted as $\mathcal{V}^-(k)$ such that $\mathcal{V}^{-}(k) = \mathcal{V} \setminus \mathcal{V}^{+}(k)$. Note that all negative valued vertices and all non-negative valued vertices that are connected to negative valued vertices are in $\mathcal{V}^{-}(k)$.
    According to the assumption provided
    at the beginning of the proof, vertices in $\mathcal{V}^+(k)$ cannot be connected to the vertices in $\mathcal{V}^-(k)$, then adjacency matrix of the time-varying weighted graph can be permuted as
    \begin{equation*}
        \widehat{A}_{a_w}(k) = P(k) A_{a_w}(k) P(k)^T = 
        \begin{bmatrix} 
            A_{a_w}^+(k) &  \bm{0} \\
            \bm{0}  &  A_{a_w}^-(k) \\
        \end{bmatrix},
    \end{equation*}
    where $P(k)$ is a permutation matrix, $A_{a_w}^+(k)$ and $A_{a_w}^-(k)$ are the adjacency matrices of sub-graphs that consist of the vertices in 
    $\mathcal{V}^{+}(k)$ and $\mathcal{V}^{-}(k)$, respectively. Note that due to Condition \ref{con:pos_w}\ref{item:c1a},
    if $A_{a_w}[i,j](k) = 0$, where $e[i](k) \geq 0$, then $A_{a}[i,j](k) = 0$. Therefore,
    $A_a$ can also be permuted as
    \begin{equation} \label{eq:reducible_adj}
        \widehat{A}_{a}(k) = P(k) A_{a} P(k)^T = 
        \begin{bmatrix} 
            A_{a}^+(k) &  \bm{0} \\
            \bm{0}  &  A_{a}^-(k) \\
        \end{bmatrix}.
    \end{equation}
    Due to Assumption \ref{asm:str_con_graph}, $\mathcal{G}$ is a connected graph, 
    the $A_a$ matrix is an irreducible matrix, and then
    $\widehat{A}_{a}(k)$ matrix is also an irreducible matrix.
    However, the matrix in Eq. (\ref{eq:reducible_adj}) is a reducible matrix, which means it contradicts the assumption provided at the beginning of the proof. Therefore, all non-negative valued vertices are connected to some negative valued vertices.
\end{proof}

\vspace{-0.12cm}
Therefore, even if $\mathcal{G}_w(k)$ may not be a connected graph for some $k$, there always exists at least one connected sub-graph that consists of both non-negative and negative valued vertices. It is important to note that in these connected sub-graphs, each edge is between a non-negative valued vertex and any other vertex and there is no edge between two negative valued vertices due to Condition \ref{con:pos_w}\ref{item:c1b}. Therefore, weights of all edges are at least $c_1$ due to Condition \ref{con:pos_w}\ref{item:c1a}.

The following lemma provides
a bound for the eigenvalues of the Laplacian matrix $\mathcal{L}_{w}(k) $ of a connected time-varying weighted graph when the weights are within a certain range.
\begin{lemma} \label{lem:gers}
    Let ${\mathcal{G}_{w}(k)}$ be a connected time-varying weighted graph with $m$ vertices. 
    Assume that the following inequality is satisfied for the weights of the edges of the graph $\mathcal{G}_{w}(k)$,
    \begin{equation} \label{eq:greg_0}
        \begin{aligned}
            &c_1 \leq w_{{i,j}}(k) \leq \frac{1-c_2}{\max(\mathcal{D} [i,i], \mathcal{D} [j,j]) }, 
        \end{aligned}
    \end{equation}
    where $0 < c_1$ and $0<c_2 < 1$. Then, the following inequality holds
    for the eigenvalues of the Laplacian matrix $\mathcal{L}_{w}(k)$,
    \begin{equation*} \label{eq:greg_3}
        \sigma( \mathcal{L}_{w}(k) ) \subseteq [c_3, (2-c_3)] \cup \{ 0 \},
    \end{equation*}
    where $0 < c_3 = \min \Big( c_1 \frac{4}{m (m-1)}, 2c_2 \Big) \leq 1$.
\end{lemma}

\begin{proof}
According to Gershgorin's Circle Theorem \cite{bell1965gershgorin}, eigenvalues of a $\mathcal{L}_{w}(k)$ matrix satisfy, $\sigma( \mathcal{L}_{w}(k) )$
\begin{equation} \label{eq:greg_1}
     \subseteq  \bigcup_i  \Big\{ \lambda \in \mathbb{C}: \big| \lambda - \mathcal{L}_{w}[i,i](k) \big| \leq \sum_{ j \in \mathcal{I}_1^m \setminus \{ i \} } \big| \mathcal{L}_{w}[i,j](k) \big| \Big \}.
\end{equation}
Since $ \mathcal{L}_{w}[i,i](k) 
= \sum_{  \{ \mathrm{v}_i, \mathrm{v}_j\} \in \mathcal{E} }  w_{{i,j}}(k) $ for weighted Laplacian matrix, where $ 0 < c_1 \leq w_{{i,j}}(k)$, Eq. (\ref{eq:greg_1}) implies that, $\sigma(  \mathcal{L}_{w}(k) ) $
\begin{equation*} 
    \begin{aligned} 
        & \subseteq \bigcup_i  \Big\{ \lambda \in \mathbb{C}: \bigg | \lambda - \sum_{ \{\mathrm{v}_i,\mathrm{v}_j\} \in \mathcal{E} }  w_{{i,j}}(k) \bigg| \leq  \sum_{ \{\mathrm{v}_i,\mathrm{v}_j\} \in \mathcal{E} }  w_{{i,j}}(k) \Big \}\\
        & = \bigcup_i  \Big\{ \lambda \in \mathbb{C}: 0 \leq \lambda \leq  2 \sum_{ \{\mathrm{v}_i,\mathrm{v}_j\} \in \mathcal{E} }  w_{{i,j}}(k) \Big \}.
    \end{aligned}
\end{equation*}
According to Eq. (\ref{eq:greg_0}), $w_{{i,j}}(k) \leq (1 - c_2) / \mathcal{D}[i,i]$ and $w_{{i,j}}(k) \leq (1 - c_2) / \mathcal{D}[j,j]$. Then, eigenvalues of the Laplacian matrix $\mathcal{L}_{w}(k)$ satisfy that,
\begin{equation} \label{eq:greg_2}
        0 \leq \lambda_i(k) \leq (2 - 2c_2) \text{ for } \lambda_i(k) \in \sigma ( \mathcal{L}_{w}(k) ).
\end{equation}

We now show that the second smallest eigenvalue of the Laplacian matrix $\mathcal{L}_{w}(k)$ is also lower bounded. 
In \cite[Theorem 4.2]{mohar1991eigenvalues}, the following lower bound is provided for the second smallest eigenvalue of the Laplacian matrix of a connected unweighted graph $\mathcal{G}$, 
\begin{equation*}
    \begin{aligned}
        \frac{4}{ \text{diam}(\mathcal{G}) m } \leq \lambda_2(\mathcal{L}), \\
    \end{aligned}
\end{equation*}
where $\text{diam}(\mathcal{G})$ is the diameter of the graph, which is the length of the longest shortest path between any two vertices, and $m$ is the number of vertices of the graph. Let us define the Laplacian matrix $\mathcal{L}_{c_1}$ for the case that weights of all edges of the graph $\mathcal{G}$ are $c_1$ instead of $1$. Suppose that, $\lambda_i(\mathcal{L})$ and $\lambda_i(\mathcal{L}_{c_1})$ are the $i^{\text{th}}$ eigenvalues of $\mathcal{L}$ and $\mathcal{L}_{c_1}$ matrices such that $\lambda_i(\mathcal{L}) \leq \lambda_{i+1}(\mathcal{L})$ and $\lambda_i(\mathcal{L}_{c_1}) \leq \lambda_{i+1}(\mathcal{L}_{c_1})$ for $i \in \mathcal{I}_1^{m-1}$. Then, $c_1 \lambda_i(\mathcal{L}) = \lambda_i(\mathcal{L}_{c_1})$ for $i \in \mathcal{I}_1^{m}$. Since $c_1 \leq w_{i,j}$ for the graph $\mathcal{G}_{w}(k)$, $\mathcal{L}_{w}(k) - \mathcal{L}_{c_1}$ is also a positive semi-definite Laplacian matrix, which means $u^T(\mathcal{L}_{w}(k) - \mathcal{L}_{c_1})u \geq 0$ for all $u \in \mathbb{R}^m$. Let us denote unit eigenvector $u_2 \in \mathbb{R}^m$ such that $\mathcal{L}_{w}(k)u_2 = \lambda_2(\mathcal{L}_{w}(k)) u_2$ and unit eigenvector $z_i \in \mathbb{R}^m$ such that $\mathcal{L}_{c_1} z_i = \lambda_i(\mathcal{L}_{c_1}) z_i$ for all $i \in \mathcal{I}_1^{m}$. Note that $\bm{1}$ is the eigenvector of these Laplacian matrices with the associated eigenvalue of $0$. Since these Laplacian matrices are symmetric, all other eigenvectors of them are orthogonal to $\bm{1}$. Then, $u_2^T z_1 = u_2^T \bm{1} = 0$ and $u_2 \in \text{span}(z_2, z_3, ..., z_m)$, which means $u_2 = \sum_{i \in \mathcal{I}_2^m}\alpha_i z_i$, where $\sum_{i \in \mathcal{I}_2^m}\alpha_i^2 = 1$. Thus,
\begin{align*}
    0 &\leq u_2^T(\mathcal{L}_{w}(k) - \mathcal{L}_{c_1})u_2\\
      &= \lambda_2(\mathcal{L}_{w}(k)) - u_2^T \mathcal{L}_{c_1} u_2\\
      &= \lambda_2(\mathcal{L}_{w}(k)) - \sum_{i \in \mathcal{I}_2^m} \alpha_i z_i^T  \lambda_i(\mathcal{L}_{c_1}) \alpha_i z_i\\
      &= \lambda_2(\mathcal{L}_{w}(k)) - \sum_{i \in \mathcal{I}_2^m} \alpha_i^2 \lambda_i(\mathcal{L}_{c_1})\\
      &\leq \lambda_2(\mathcal{L}_{w}(k)) - \lambda_2(\mathcal{L}_{c_1}).
\end{align*}
Hence, Eq. (\ref{eq:Lap_w}) can be provided for the second smallest eigenvalue of $\mathcal{L}_{w}(k)$,
\begin{equation} \label{eq:Lap_w}
    c_1 \begin{aligned}
        \frac{4}{ \text{diam}(\mathcal{G}) m } \leq \lambda_2(\mathcal{L}_{w}(k)).
    \end{aligned}
\end{equation}
Since $0 < \text{diam}(\mathcal{G}) \leq m-1$ for a connected graph and due to Eq. (\ref{eq:greg_2}), the following equation represents the eigenvalues of the matrix $ \mathcal{L}_{w}(k)$,
\begin{equation*}
    \sigma( \mathcal{L}_{w}(k)) \subseteq [c_3, (2-c_3)] \cup \{ 0 \},
\end{equation*}
where $0 < c_3 = \min \Big( c_1 \frac{4}{m (m-1)}, 2c_2 \Big) \leq \min ( 2c_1 , 2c_2 ) \leq \min ( 2 - 2c_2 , 2c_2 ) < 1$. 
\end{proof}
\vspace{0.3cm}

The following theorem shows that values of vertices of $\mathcal{G}_{w}(k)$ converge to $0$ under specific value transition dynamics.
\begin{theorem} \label{thm:graph_general_2}
Suppose that Assumption \ref{asm:str_con_graph} and Condition \ref{con:pos_w} are satisfied, and $e(k)$ vector in Definition \ref{def:graph_def} evolves according to 
\begin{equation} \label{eq:trans_graph}
    \begin{aligned} 
        e[i](k+1) = e[i](k) + \sum_{j \in \mathcal{I}_1^m } \Big( e[j](k) - e[i](k) \Big) A_{a_w} [i,j] (k), 
    \end{aligned}
\end{equation}
for all $i \in \mathcal{I}_1^m$, where 
$A_{a_w}[i,j](k) = A_a[i,j]  w_{i,j}(k), \;
0 \leq w_{i,j}(k) \leq  (1 - c_2)/(\max(\mathcal{D} [i,i], \mathcal{D} [j,j])), \; 0 < c_2 < 1$. Then $e(k)$ exponentially convergences to $\bm{0}$, i.e. $ ||e(k+1)|| / ||e(k)|| \leq \gamma$, where $0 \leq \gamma < 1 $.
\end{theorem}

\begin{proof}
Eq. (\ref{eq:trans_graph}) is equivalent to following equation,
\begin{equation*}
    \begin{aligned} \label{eq:Lin_Discrite}
        e(k+1) = F&(k) e(k), \text{ where} \\ 
        F(k) = I -& \mathcal{L}_w(k). \\
    \end{aligned}
\end{equation*}

Note that $F(k)$ is a doubly stochastic symmetric matrix, where $\bm{1}^T F(k) = \bm{1}^T$, $F(k)\bm{1} = \bm{1}$, $F(k)^T = F(k)$ and non-zero, non-diagonal elements of $F(k)$ represents the weights of activated edges. 

Condition \ref{con:pos_w}\ref{item:c1a} forces weights of the edges connecting non-negative valued vertices to other vertices to be at least $c_1$. However, the weights of some other edges
may be $0$. Thus, the weighted graph may not be connected, which means $\mathcal{L}_w(k)$ and $F(k)$ matrices may be reducible. Hence, $\mathcal{L}_w(k)$ may have repeated eigenvalues at $0$, and $F(k)$ may have repeated eigenvalues at $1$. 

Even if $\mathcal{G}_w(k)$ may not be connected, due to Lemma \ref{lem:err_vec}, all non-negative valued vertices belong to a connected sub-graph in $\mathcal{G}_w(k)$, which consists of both non-negative and negative valued vertices.
Note that in these connected sub-graphs,
there is no edge between two negative valued vertices due to Condition \ref{con:pos_w}\ref{item:c1b}. Therefore, weights of all edges
are at least $c_1$ due to Condition \ref{con:pos_w}\ref{item:c1a}, which enables us to apply Lemma \ref{lem:gers} in later parts.
Let us denote the number of these connected sub-graphs as $s(k) \geq 1$ and denote the reordered $F(k)$ matrix and $e(k)$ vector as
\begin{equation*}
    \begin{aligned}
    \widehat{F}(k) &= P(k) F(k) P(k)^T \\
    &=  \begin{bmatrix} 
        \widehat{F}_{0}(k)  &                           &                       &        \\
                            &\widehat{F}_{1}(k)&        &                                \\
                            &                           & \ddots &                       \\
                            &                           &        & \widehat{F}_{s(k)}(k) 

        \end{bmatrix}, \\
    \hat{e}(k) &= P(k) e(k) \\
               &= (\hat{e}_{0}(k),\hat{e}_1(k),..., \hat{e}_{s(k)}(k)),
    \end{aligned}
\end{equation*}
where $P(k)$ is a permutation matrix, $\widehat{F}_i(k)$, $i \in I_s(k) = \{ 1,...,s(k) \}$ are irreducible matrices that correspond to connected sub-graphs and $\widehat{F}_0(k)$ is a matrix that corresponds to remaining sub-graph.
Note that there are both non-negative and negative elements in $\hat{e}_i(k)$, $i \in I_s(k)$ and all elements of $\hat{e}_0(k)$ are negative.

Recall that, in Example 1, $s(k) = 1$, values of $\mathrm{v}_1$, $\mathrm{v}_2$ and $\mathrm{v}_4$ vertices belong to the $\hat{e}_1(k)$ vector and values of $\mathrm{v}_3$, $\mathrm{v}_5$ and $\mathrm{v}_6$ vertices belong to the $\hat{e}_0(k)$ vector.

Suppose that, $\alpha_i(k)$ is the average of elements in $\hat{e}_i(k)$ such that $\alpha_i(k) = \bm{1}^T \hat{e}_i(k)/m_i(k)$, where $m_i(k)$ is the number of elements in $\hat{e}_i(k)$ and $\bm{\alpha}_i(k) =  \bm{1} \alpha_i(k) $, $\bm{\alpha}_i(k) \in \mathbb{R}^{m_i(k)}$. Lemma \ref{lem:gers} implies that only one eigenvalue of $\widehat{F}_i(k)$, $i \in I_s(k)$ is $1$ with associated eigenvector $\bm{1}$, and other eigenvalues are in $[-(1-c_{3_i}), (1-c_{3_i})]$, $0 < c_{3_i} \leq 1$ since $\widehat{F}_i(k)$, $i \in I_s(k)$ matrices corresponds to a connected sub-graphs in which weights of all edges are at least $c_1 > 0$. Since the values of all vertices are negative in the remaining sub-graph that corresponds to the $\widehat{F}_0(k)$ matrix, there is no lower-bound for the weights of the edges, while there is an upper bound provided in Eq. (\ref{eq:trans_graph}). Hence, $\sigma(\widehat{F}_0(k)) \in [-(1-c_{3_0}),1]$, $0 < c_{3_0} \leq 1$. Therefore the following equation holds,
\begin{equation} \label{eq:dec_3}
    \begin{aligned} 
        & ||\widehat{F}_i(k) \hat{e}_i(k) - \bm{\alpha}_i(k) ||^2 = \\
        & \qquad =||\widehat{F}_i(k) \hat{e}_i(k) - \widehat{F}_i(k) \bm{\alpha}_i(k) ||^2\\
        & \qquad =||\widehat{F}_i(k) (\hat{e}_i(k) - \bm{\alpha}_i(k)) ||^2\\
        & \qquad = (\hat{e}_i(k) - \bm{\alpha}_i(k))^T \widehat{F}_i(k)^T \widehat{F}_i(k) (\hat{e}_i(k) - \bm{\alpha}_i(k))\\
        & \qquad \leq c_4^2 (\hat{e}_i(k) - \bm{\alpha}_i(k))^T (\hat{e}_i(k) - \bm{\alpha}_i(k))\\
        & \qquad = c_4^2 ||\hat{e}_i(k) - \bm{\alpha}_i(k) ||^2,
    \end{aligned}
\end{equation}
where $c_4^2 = (1- \min_i c_{3_i} )^2 < 1$ for all $i \in I_s(k)$.

Here, we will show that $||\widehat{F}_i(k) \hat{e}_i(k)||^2$ is also smaller than or equal to $c_5^2||\hat{e}_i(k)||^2$ for all $i \in \mathcal{I}_1^m$, where $c_5^2<1$. Note that if $ \bm{1}^T \bm{\alpha}_i(k) = \bm{1}^T \hat{e}_i(k) $, then $\bm{\alpha}_i(k)$ and  $(\hat{e}_i(k) - \bm{\alpha}_i(k))$ are orthogonal to each other, 
\begin{equation*}
    \begin{aligned}
        \bm{\alpha}_i(k)^T (\hat{e}_i(k) - \bm{\alpha}_i(k)) &= \bm{\alpha}_i(k)^T \hat{e}_i(k) - \bm{\alpha}_i(k)^T \bm{\alpha}_i(k)\\
        &=\alpha_i(k) \bm{1}^T \hat{e}_i(k) - \bm{\alpha}_i(k)^T \bm{\alpha}_i(k)\\
        &= \alpha_i(k) \bm{1}^T \bm{\alpha}_i(k) - \bm{\alpha}_i(k)^T \bm{\alpha}_i(k)\\
        &= \bm{\alpha}_i(k)^T \bm{\alpha}_i(k) - \bm{\alpha}_i(k)^T \bm{\alpha}_i(k)\\
        &= 0.
    \end{aligned}
\end{equation*}

Hence, it can be shown that,
\begin{equation*}
    \begin{aligned}
        ||\hat{e}_i(k) - \bm{\alpha}_i(k) ||^2 &= ||\hat{e}_i(k)||^2 -2\bm{\alpha}_i(k)^T \hat{e}_i(k) + ||\bm{\alpha}_i(k)||^2\\
        &=||\hat{e}_i(k)||^2 -2\bm{\alpha}_i(k)^T \bm{\alpha}_i(k) + ||\bm{\alpha}_i(k)||^2\\
        &=||\hat{e}_i(k)||^2 - ||\bm{\alpha}_i(k)||^2,
    \end{aligned}
\end{equation*}
and $||\widehat{F}_i(k)\hat{e}_i(k) - \bm{\alpha}_i(k) ||^2 = $
\begin{equation*}
    \begin{aligned}
        &= ||\widehat{F}_i(k)\hat{e}_i(k)||^2 -2\bm{\alpha}_i(k)^T \widehat{F}_i(k)\hat{e}_i(k) + ||\bm{\alpha}_i(k)||^2\\
        &= ||\widehat{F}_i(k)\hat{e}_i(k)||^2 -2 \alpha_i(k) \bm{1}^T \widehat{F}_i(k)\hat{e}_i(k) + ||\bm{\alpha}_i(k)||^2\\
        &= ||\widehat{F}_i(k)\hat{e}_i(k)||^2 -2 \alpha_i(k) \bm{1}^T \hat{e}_i(k) + ||\bm{\alpha}_i(k)||^2\\
        &= ||\widehat{F}_i(k)\hat{e}_i(k)||^2 -2 \alpha_i(k) \bm{1}^T \bm{\alpha}_i(k)^T + ||\bm{\alpha}_i(k)||^2\\
        &=||\widehat{F}_i(k)\hat{e}_i(k)||^2 -2\bm{\alpha}_i(k)^T \bm{\alpha}_i(k) + ||\bm{\alpha}_i(k)||^2\\
        &=||\widehat{F}_i(k)\hat{e}_i(k)||^2 - ||\bm{\alpha}_i(k)||^2.
    \end{aligned}
\end{equation*}

If $i \in I_s(k)$, then $\hat{e}_i(k) \neq \bm{0}$ because there are both non-negative and negative values in $\hat{e}_i(k)$.
Therefore, Eq. (\ref{eq:dec_3}) implies that
\begin{equation} \label{eq:dec_4}
    \begin{aligned}
        ||\widehat{F}_i(k) \hat{e}_i(k) - \bm{\alpha}_i(k) ||^2 &\leq c_4^2 ||\hat{e}_i(k) - \bm{\alpha}_i(k) ||^2\\
        ||\widehat{F}_i(k) \hat{e}_i(k)||^2 - ||\bm{\alpha}_i(k) ||^2 &\leq c_4^2 (||\hat{e}_i(k)||^2 - ||\bm{\alpha}_i(k) ||^2)\\
        ||\widehat{F}_i(k) \hat{e}_i(k)||^2 &\leq c_4^2 ||\hat{e}_i(k)||^2 + (1-c_4^2) ||\bm{\alpha}_i(k) ||^2\\
        \frac{||\widehat{F}_i(k) \hat{e}_i(k)||^2}{||\hat{e}_i(k)||^2} &\leq c_4^2 + (1-c_4^2) \frac{||\bm{\alpha}_i(k) ||^2}{||\hat{e}_i(k)||^2}.\\
    \end{aligned}
\end{equation}

Here, we will show that $\frac{||\bm{\alpha}_i(k) ||^2}{||\hat{e}_i(k)||^2} < \frac{m-1}{m}$ for all $i \in I_s(k)$. Suppose that, $\hat{e}_i(k) = P_{\hat{e}_i}(k) (\hat{e}_{i_p}(k), \hat{e}_{i_n}(k))$, where $P_{\hat{e}_i}(k)$ is a permutation matrix such that $\hat{e}_{i_p}(k) \geq 0$ and $\hat{e}_{i_n}(k) < 0$. Let $m_{i_p}(k)$ and $m_{i_n}(k)$ represent the number of elements in $\hat{e}_{i_p}(k)$ and $\hat{e}_{i_n}(k)$, where $m_{i_p}(k) + m_{i_n}(k) = m_{i}(k)$. Since there are both non-negative and negative values in $\hat{e}_{i}(k)$, $m_{i_p}(k) \geq 1$, $m_{i_n}(k) \geq 1$ and $\bm{1}^T \hat{e}_i(k) < ||\hat{e}_i(k)||_1$, then
\begin{align*}
        \frac{||\bm{\alpha}_i(k) ||^2}{||\hat{e}_i(k)||^2} &= \frac{ \bm{\alpha}_i^T \bm{\alpha}_i }{||\hat{e}_i(k)||^2}\\
        &= \frac{ \bm{1}^T \bm{1} (\bm{1}^T e_i(k))^2 }{ m_i^2(k) ||\hat{e}_i(k)||^2}\\
        &= \frac{ (\bm{1}^T e_i(k))^2 }{m_i(k) ||\hat{e}_i(k)||^2}\\
        &< \frac{ (\bm{1}^T \hat{e}_{i_p}(k))^2 + (\bm{1}^T \hat{e}_{i_n}(k))^2  }{ m_i(k) ||\hat{e}_i(k)||^2}\\
        &= \frac{ ||\hat{e}_{i_p}(k)||_1^2 + ||\hat{e}_{i_n}(k)||_1^2}{ m_i(k) ||\hat{e}_i(k)||^2}\\
        &\leq \frac{ m_{i_p}(k) ||\hat{e}_{i_p}(k)||^2 + m_{i_n}(k) ||\hat{e}_{i_n}(k)||^2}{ m_i(k) ||\hat{e}_i(k)||^2}\\
        &\leq \frac{ \max(m_{i_p}(k), m_{i_n}(k))  (||\hat{e}_{i_p}(k)||^2 + ||\hat{e}_{i_n}(k)||^2) }{ m_i(k) ||\hat{e}_i(k)||^2}\\
        &= \frac{ \max(m_{i_p}(k), m_{i_n}(k)) }{ m_i(k) }\\
        &\leq \frac{ m_i(k) - 1  }{ m_i(k) } \leq \frac{ m - 1  }{ m }\\
\end{align*}
for all $i \in I_s(k)$. It then follows that Eq. (\ref{eq:dec_4}) implies
\begin{equation*}
    \begin{aligned}
        ||\widehat{F}_i(k) \hat{e}_i(k)||^2 &\leq c_5^2 ||\hat{e}_i(k)||^2 \text{\; for all \;} i \in I_s(k),\\
    \end{aligned}
\end{equation*}
where $c_5^2 < c_4^2 + (1 - c_4^2) \frac{ m - 1  }{ m } < 1 $.

Here, we will present that $||e(k+1)|| \leq c_6||e(k)||$, where $0<c_6<1$. Since $e(k) \neq \bm{0}$,
\begin{equation*}
    \begin{aligned}
             \frac{||e(k+1)||^2}{||e(k)||^2} 
              =  \frac{||\widehat{F}(k)\hat{e}(k)||^2 }{||\hat{e}(k)||^2}   = \frac{ \sum_{i=0}^{s(k)} ||\widehat{F}_i(k)\hat{e}_i(k)||^2 }{ \sum_{i=0}^{s(k)} ||\hat{e}_i(k)||^2}.
    \end{aligned}
\end{equation*}

Since $\sigma(\widehat{F}_0(k)) \subseteq [-(1-c_{3_0}),1]$, $0 < c_{3_0} \leq 1$, the equation is derived as
\begin{equation} \label{eq:conv_mid}
    \begin{aligned}
             \frac{||e(k+1)||^2}{||e(k)||^2} \leq \frac{ || \hat{e}_0(k)||^2 + c_5^2 \sum_{i=1}^{s(k)} ||\hat{e}_i(k)||^2 }{ ||\hat{e}_0(k)||^2 + \sum_{i=1}^{s(k)} ||\hat{e}_i(k)||^2}.
    \end{aligned}
\end{equation}

Here, we derive an upper bound for $||\hat{e}_0(k)||^2$ with respect to other terms to show that,
there is an upper bound for the right side of Eq. (\ref{eq:conv_mid}).
Since $\bm{1}^T e(k) = 0$ and $\hat{e}_0(k) < \bm{0}$, $|| \hat{e}_{0}(k) ||_1 + \sum_{i=1}^{s(k)} || \hat{e}_{i_n}(k) ||_1$  =  $\sum_{i=1}^{s(k)} || \hat{e}_{i_p}(k) ||_1$. Note that $|| \hat{e}_{i_n}(k) ||_1 > 0$ for any $i \in I_s(k)$ due to Lemma \ref{lem:err_vec}. Then,

\begin{equation} \label{eq:e_1_2}
    \begin{aligned} 
        ||\hat{e}_0(k)||_1 &< \sum_{i=1}^{s(k)} ||\hat{e}_{i_p}(k)||_1 + \sum_{i=1}^{s(k)} ||\hat{e}_{i_n}(k)||_1\\ 
        &= \sum_{i=1}^{s(k)} ||\hat{e}_i(k)||_1.
    \end{aligned}
\end{equation}

The following inequalities can be provided using equivalence of norms,
\begin{equation*} \label{eq:eq_norms}
    \begin{aligned}
        ||\hat{e}_0(k)||^2 \leq ||\hat{e}_0(&k)||_1^2 \leq m_0 ||\hat{e}_0(k)||^2, \\
        \sum_{i=1}^{s(k)}||\hat{e}_i(k)||^2 \leq \sum_{i=1}^{s(k)}||\hat{e}_i(&k)||_1^2 \leq (m-m_0) \sum_{i=1}^{s(k)}||\hat{e}_i(k)||^2.
    \end{aligned}
\end{equation*}

Thus, Eq. (\ref{eq:e_1_2}) implies that
\begin{equation*}
    ||\hat{e}_0(k)||^2 < (m-m_0)\sum_{i=1}^{s(k)} ||\hat{e}_i(k)||^2.
\end{equation*}

The strict upper-bound for $||\hat{e}_0(k)||^2$ can be inserted into Eq. (\ref{eq:conv_mid}) as
\begin{equation*}
    \begin{aligned}
        \frac{||e(k+1)||^2}{||e(k)||^2} &\leq 1 - \frac{ (1 - c_5^2) \sum_{i=1}^{s(k)} ||\hat{e}_i(k)||^2 }{ ||\hat{e}_0(k)||^2 + \sum_{i=1}^{s(k)} ||\hat{e}_i(k)||^2}\\
        &<  1 - \frac{ (1 - c_5^2) \sum_{i=1}^{s(k)} ||\hat{e}_i(k)||^2 }{ (m-m_0+1) \sum_{i=1}^{s(k)} ||\hat{e}_i(k)||^2}\\
        &= \frac{ m - m_0 + c_5^2 }{ m - m_0 + 1 }.
    \end{aligned}    
\end{equation*}

Since $0 \leq c_5^2 < 1$ and $0 \leq m_0 \leq m - 2$,
\begin{equation*}
    \frac{||e(k+1)||}{||e(k)||} \leq \gamma, 
\end{equation*}
where $\gamma < \sqrt{\max_{0 \leq m_0 \leq m-2} \frac{ m - m_0 + c_5^2 }{ m - m_0 + 1 }} = \sqrt{\frac{m + c_5^2}{m+1}} < 1$.
Therefore, $||e(k)||$ exponentially converges to $0$. 
\end{proof}
\section{Decentralized State-Dependent \\ Markov Chain Synthesis} \label{sec:somm}
Based on the consensus protocol developed in Section \ref{sec:DSDN},
we propose the decentralized state-dependent Markov chain synthesis (DSMC) algorithm that achieves convergence to the desired distribution with an exponential rate and minimal state transitions.
Additionally, we present a shortest path algorithm that can be integrated with the DSMC algorithm, as utilized in \cite{accikmecse2012markov, bandyopadhyay2017probabilistic, jang2018local}, to further enhance the convergence rate.

\vspace{-0.25cm}
\subsection{The DSMC Algorithm} \label{sec:dsmc_1}
We define a finite-state discrete-time Markov chain evolving on the vertices of the uniform graph $\mathcal{G}$ in Definition \ref{def:graph_def}.
\vspace{-0.12cm}
\begin{definition} \label{def:Markov}
    \textit{(Markov chain)}
    Let $\mathcal{V} = \{ \mathrm{v}_1, ..., \mathrm{v}_m \}$, which is the vertices of the graph in Definition \ref{def:graph_def}, be the finite set of states of the Markov chain and $\mathcal{X}(k) \in \mathcal{V}$ be the random variable of the Markov chain for time $k \in \mathbb{Z}_+$.
    The connectivity of the states is determined by the adjacency matrix $A_a$ of the uniform graph $\mathcal{G}$. Additionally, we define:
    \begin{itemize}
        \item (\textit{Probability distribution})
        Let $x(k) \in \mathbb{R}^m$ be the probability distribution at time $k: x[i](k) = \mathcal{P}(\mathcal{X}(k) = \mathrm{v}_i)$. 
        \item (\textit{Markov matrix})
           Markov matrix determines the state transitions of the Markov chain $M(k) \in \mathbb{R}^{m \times m}$,
            where $M[i,j](k) = \mathcal{P}(\mathcal{X}(k+1) = \mathrm{v}_i | \mathcal{X}(k) = \mathrm{v}_j )$ for all $i, j \in \mathcal{I}_1^m$. 
            Hence $x(k+1) = M(k)x(k)$ for $k \in \mathbb{Z}_+$.
            \item (\textit{Desired steady-state distribution})
        There exists a prescribed  desired steady-state probability distribution  $v \in \mathbb{R}^m$ such that  
        \begin{equation*} 
            \lim_{k\to\infty} x(k) = v.
        \end{equation*} 
    \end{itemize}
\end{definition}
\vspace{-0.12cm}

Let the error vector $e(k)$ be the difference between the probability distribution at time $k$ and the desired steady-state distribution  $e(k) = x(k) - v$.
The DSMC algorithm is designed to ensure that the dynamics of the error vector are identical to the dynamics of the value vector described in Theorem \ref{thm:graph_general_2}. This design guarantees consistency between the two, ensuring desirable convergence properties and performance in the DSMC algorithm.
Similar to the consensus protocol in Section \ref{sec:DSDN}, transitions in the DSMC algorithm occur from the states with higher error values to their adjacent states with lower error values to equalize the error values across all states of the Markov chain. Since the sum of these error values is equal to $0$,  the error values at all states will eventually become balanced at $0$,
resulting in the convergence of the probability distribution to the desired distribution.

\begin{algorithm}[!hbt]
    \textbf{Input:} $x[l](k), v[l], \, \forall l \in I_{ij}$,
    $A_{a}[\ell, \cdot ]$,
    $\forall \ell \in \{ i, j \} $, $c_2 \in (0,1)$\\
    \textbf{Output:} $M[i,j](k)$
    \begin{algorithmic}[1]
        \State Compute error values at the states associated with set $I_{ij}$:
        
        $$e[l](k) := x[l](k) - v[l], \, \forall l = I_{ij}$$
        \State Calculate the non-negative error differences between adjacent states:
        $$E[i,j](k) := \max \bigg( 0, \Big( e[j](k) - e[i](k) \Big)  A_{a}[i,j] \bigg)$$
        \State Determine the total number of adjacent states:
        $$d[\ell] := \sum_{\substack{ l \in \mathcal{I}_1^{m} \setminus \ell }} A_{a}[\ell,l], \, \forall \ell \in \{ i,j \}$$
        \State Calculate the scaling factor to ensure algorithm stability:
        $$s_{1}[i,j] := \frac{1-c_2}{\max(d[i], d[j])}$$
        \State Set the scaling factor based on the probability of state $\mathrm{v}_j$:
        $$s_{2}[j](k) := 
            \begin{cases}
                \frac{x[j](k) }{ \sum_{l \in \mathcal{I}_1^{m} } E[l,j](k)} & \text{if } \sum_{l \in \mathcal{I}_1^{m} } E[l,j](k) \neq 0 \\
                1                                           & \text{otherwise}
            \end{cases} $$
        \State Determine the scaling factor to satisfy Condition \ref{con:pos_w}\ref{item:c1b}:
        $$s_{3}[i,j](k) := 
                \begin{cases}
                    0                                           & 
                    
                    \text{\parbox{15em}{if $e[\ell] < 0$ for all $\ell \in \{ i,j \}$ and $e[l](k) \geq 0$ for some $l \in I_{ij}$}}\\
                    1                                           & \text{otherwise}
                \end{cases}$$
        \State Obtain the resulting scaling factor:
        $$W[i,j](k) := \min\Big(s_1[i,j] , s_2[j](k) , s_3[i,j](k) \Big)$$
        \State Calculate the scaled non-negative error differences between adjacent states:
        $$T[i,j](k) := E[i,j](k) W[i,j](k)$$
        \State Compute the transition probability from $\mathrm{v}_j$ to $\mathrm{v}_i$:
        $$F[i,j](k) := 
            \begin{cases}
                \frac{T[i,j](k)}{x[j](k)} & \text{if $x[j](k) \neq 0$} \\
                0 & \text{otherwise}
            \end{cases}$$
        \State Synthesize the Markov matrix:
        $$M[i,j](k) :=
            \begin{cases}
                F[i,j](k) & \text{if $i \neq j$} \\
                1 - \sum_{l \in \mathcal{I}_1^{m} } F[l,j](k) & \text{if $i = j$}
            \end{cases}$$
    \end{algorithmic}
    \caption{Decentralized state-dependent Markov chain synthesis algorithm}
    \label{alg:DSMC}
\end{algorithm}

We present the synthesis of the Markov matrix in Algorithm \ref{alg:DSMC}, and its convergence proof is presented in Theorem \ref{thm:dsmc_conv}.
To this end, let $I_{ij}$ be the set of indices of states adjacent to either $\mathrm{v}_i$ or $\mathrm{v}_j$, such that $I_{ij} = \{ l : A_a[l,i] = 1 \text{ or } A_a[l,j] = 1 \}$.
It should be noted that $e[\ell] < 0$ for all $\ell \in \{ i, j \}$ and $e[l](k) \geq 0$ for some $l \in I_{ij}$ if and only if $ \{i, j\} \in I_{n_p, n}(k)$, which is introduced in Definition \ref{def:erridxset}.
To determine the transition probability from state $\mathrm{v}_j$ to state $\mathrm{v}_i$ for time $k$, the required inputs of the algorithm are $x[l](k)$ and $v[l]$ for all $l \in I_{ij}$, the corresponding rows of the adjacency matrix for states $\mathrm{v}_i$ and $\mathrm{v}_j$, and a hyper-parameter $c_2 \in (0,1)$ that scales the transition probabilities.
In Line 1 of Algorithm \ref{alg:DSMC}, 
the algorithm computes the error values at the states
associated with the set $I_{ij}$.
From the state $\mathrm{v}_j$, there will be a transition to any adjacent state $\mathrm{v}_i$ only if $e[j](k) > e[i](k)$.
Also, the transition probabilities
vary proportionally with the differences in error values.
Therefore, the difference of the error value $e[j](k)$ from the error value $e[i](k)$ is calculated in Line 2 if $e[j](k) > e[i](k)$, otherwise it is set to $0$.
As in Theorem \ref{thm:graph_general_2}, appropriately scaling the differences between error values relative to the number of adjacent states is crucial for the stability of the convergence.
Also, it is required to scale these differences with respect to the probabilities of the corresponding states.
For this purpose, a scaling factor is determined in Line 4 for the stability of the convergence,
while another scaling factor is determined by the probability of the state $\mathrm{v}_j$ in Line 5.
Furthermore, if $\sum_{l \in \mathcal{I}_1^{m} } E[l,j](k) = 0$, then there will be no transition from state $\mathrm{v}_j$ to any adjacent state $\mathrm{v}_i$. In this case, the value of $s_2[j](k)$ becomes irrelevant and is set to $1$.
If $e[i](k) < 0$ and $e[j](k) < 0$, and either $\mathrm{v}_i$ or $\mathrm{v}_j$ has an adjacent state with a positive error value, then another scaling factor $s_3[i,j](k)$ is set to $0$ in Line 6 to satisfy Condition \ref{con:pos_w}\ref{item:c1b}.
The minimum of these three scaling factors is chosen in Line 7
to satisfy all conditions simultaneously. 
Hence, the scaled difference of the error value $e[j](k)$ from the error value $e[i](k)$ is calculated in Line 8.
Note that $T[i,j](k)$ represents the amount of increase in the error value at state $\mathrm{v}_i$ due to the transition from state $\mathrm{v}_j$. Equivalently, it represents the amount of decrease of the error value at state $\mathrm{v}_j$ due to the transition to the state $\mathrm{v}_i$.
This amount is divided by the probability of the state $\mathrm{v}_j$ in Line 9
to decide the transition probability from state $\mathrm{v}_j$ to state $\mathrm{v}_i$.
In Line 10,
the remaining probabilities are added to the diagonal elements, finalizing the synthesis of the Markov matrix. The complexity of the DSMC algorithm is $O(m^2)$, where $m$ is the number of states.

The following proposition shows that the matrix synthesized by Algorithm \ref{alg:DSMC} is a column stochastic Markov matrix.
\vspace{-0.20cm}
\begin{proposition} \label{pro:markov_prop}
The matrix synthesized by Algorithm \ref{alg:DSMC} satisfies $M(k) \geq \bm{0}$ and $\bm{1}^T M(k) = \bm{1}^T$ constraints.
\end{proposition}

\begin{proof}
    (Proof of $M(k) \geq \bm{0}$) $W(k) \geq \bm{0}$ since $s_1[i,j] > 0$, $s_2[j](k) \geq 0$ and $s_3[i,j](k) \geq 0$, $\forall i,j \in \mathcal{I}_1^m$. Since $E(k) \geq \bm{0}$ and $W(k) \geq \bm{0}$, $T(k) \geq \bm{0}$. $x(k) \geq \bm{0}$ and $T(k) \geq \bm{0}$ imply that $M[i,j](k) = F[i,j](k) \geq 0$, $ \forall i,j \in \mathcal{I}_1^m$ if $i \neq j$.
    
    If $x[j](k) = 0$ or $\sum_{i \in \mathcal{I}_1^{m}} E[i,j](k) = 0$, then $\sum_{i \in \mathcal{I}_1^{m}} F[i,j](k) = 0$ and $M[j,j](k) = 1$. If $x[j](k) \neq 0$ and $\sum_{i \in \mathcal{I}_1^{m}} E[i,j](k) \neq 0$, then
    \begin{equation*}
        \begin{aligned}
            \sum_{i \in \mathcal{I}_1^{m}} F[i,j](k) &= \frac{\sum_{i \in \mathcal{I}_1^{m}} T[i,j](k)}{x[j](k)}\\
            &= \frac{\sum_{i \in \mathcal{I}_1^{m}} E[i,j](k)W[i,j](k)}{x[j](k)}\\
            &\leq  \frac{\sum_{i \in \mathcal{I}_1^{m}} E[i,j](k)}{x[j](k)} s_2[j](k)\\
            &= \frac{\sum_{i \in \mathcal{I}_1^{m}} E[i,j](k)}{x[j](k)} \frac{x[j](k)}{\sum_{i \in \mathcal{I}_1^{m}} E[i,j](k)}\\
            &=1.
        \end{aligned}
    \end{equation*}
    Therefore, $\sum_{i \in \mathcal{I}_1^{m}} F[i,j](k) \leq 1$ and $M[j,j](k) \geq 0$.

    (Proof of $\bm{1}^T M(k) = \bm{1}^T$) $\bm{1}^T M(k) = \bm{1}^T$ requires that $\sum_{i \in \mathcal{I}_1^{m}} M[i,j](k) = 1$, $\forall j \in \mathcal{I}_1^{m}$. Note that $\sum_{ i \in \mathcal{I}_1^{m} \setminus \{ j \} } F[i,j](k) = \sum_{i \in \mathcal{I}_1^{m}} F[i,j](k)$ since $F[i,i](k) = 0$. Then, $\sum_{i \in \mathcal{I}_1^{m}} M[i,j](k) = \sum_{i \in \mathcal{I}_1^{m} \setminus \{ j \}} F[i,j](k) + 1 - \sum_{ i \in \mathcal{I}_1^{m} } F[i,j](k) = 1$. This shows that $\bm{1}^T M(k) = \bm{1}^T$.
\end{proof}

The following proposition shows that the error dynamics 
in Algorithm \ref{alg:DSMC} are identical
to the evolution of the value vector in the consensus protocol presented in Theorem \ref{thm:graph_general_2}. 
\vspace{-0.20cm}
\begin{proposition} \label{pro:error_evol}
    The error dynamics of Algorithm \ref{alg:DSMC} satisfy Eq. (\ref{eq:trans_graph}).
\end{proposition}
\begin{proof}
    In Algorithm \ref{alg:DSMC}, the increase in the error value at state $\mathrm{v}_i$ caused by state $\mathrm{v}_j$ is denoted by $T[i,j](k)$ in Line 8. Equivalently, $T[i,j](k)$ represents the decrease in the error value at state $\mathrm{v}_j$ caused by state $\mathrm{v}_i$.
    Since $\sum_{ i \in \mathcal{I}_1^{m} } F[i,j](k) \leq 1$ in Line 9 of Algorithm \ref{alg:DSMC}, the amount of change determined in Line 8 is not suppressed by $x(k)$. Thus, the error value at the state $\mathrm{v}_i$ changes as
    \begin{equation} \label{eq:trans_alg_1}
        \begin{aligned}
            e[i](k+1) &= e[i](k) + \sum_{j \in \mathcal{I}_1^{m}} T[i,j](k) - \sum_{j \in \mathcal{I}_1^{m}} T[j,i](k)\\
            = e[i](k) +& \hspace{-0.15cm}\sum_{j \in \mathcal{I}_1^{m}} \hspace{-0.1cm} E[i,j](k)W[i,j](k) - \hspace{-0.15cm}\sum_{j \in \mathcal{I}_1^{m}} \hspace{-0.1cm} E[j,i](k)W[j,i](k) \\
            = e[i](k) +& \sum_{\substack{ j \in \mathcal{I}_1^{m} \\ e[j](k) \geq e[i](k)}} \Big(e[j](k) - e[i](k)\Big)A_{a}[i,j] W[i,j](k) \\
            -& \sum_{\substack{ j \in \mathcal{I}_1^{m} \\ e[j](k) < e[i](k)}} \Big(e[i](k) - e[j](k)\Big)A_{a}[j,i] W[j,i](k).
        \end{aligned}
    \end{equation}
    Let us denote $w_{i,j}(k)$ as 
    \begin{equation} \label{eq:w_to_W}
        w_{i,j}(k) = W[i,j](k) \text{ if } e[j](k) \geq e[i](k).
    \end{equation}
    Note that $0 \leq W[i,j](k) \leq s_1[i,j] =(1-c_2) / \max(d[i], d[j])$, for all $i,j \in \mathcal{I}_1^m$, where $0 < c_2 < 1$ due to Line 4 and Line 7 of the Algorithm \ref{alg:DSMC}.
    Since $d[i]$ in Line 3 of the Algorithm \ref{alg:DSMC} represents the degree similar to $\mathcal{D}[i,i]$ in Definition \ref{def:graph_def}, $w_{i,j}(k)$ satisfies the constraints given in Theorem \ref{thm:graph_general_2}. Thus, Eq. (\ref{eq:trans_alg_1}) implies that
    \begin{equation} \label{eq:trans_alg_2}
        \begin{aligned}
            e[i](k+1) = e[i](k) & + \sum_{ j \in \mathcal{I}_1^{m} } \Big( e[j](k) - e[i](k) \Big) A_{{a}_w} [i,j] (k) 
        \end{aligned}
    \end{equation}
    for all $i \in \mathcal{I}_1^m$, where  $A_{{a}_w}[i,j](k) = A_{a}[i,j]  w_{i,j}(k), \; 0 \leq w_{i,j}(k) < (1 - c_2) / (\max(\mathcal{D} [i,i], \mathcal{D} [j,j])), \; 0 < c_2 < 1$.
    Since Eq. (\ref{eq:trans_alg_2}) is identical to Eq. (\ref{eq:trans_graph}), the error dynamics of Algorithm \ref{alg:DSMC} are identical to the evolution of the value vector in Theorem \ref{thm:graph_general_2}.
\end{proof}

Note that Assumption \ref{asm:str_con_graph} is already satisfied for the DSMC algorithm due to Assumption \ref{asm:strongly_connected}. We now show that Condition \ref{con:pos_w}, which is crucial for the convergence of the proposed consensus protocol, is satisfied by Algorithm \ref{alg:DSMC}.
\begin{proposition} \label{pro:ams_satis}
    Condition \ref{con:pos_w}
    is satisfied by Algorithm \ref{alg:DSMC}.
\end{proposition}
\begin{proof}
    To ensure Condition \ref{con:pos_w}\ref{item:c1a} it must be shown that there exists a $c_1$ such that $0 < c_1 \leq w_{i,j}(k)$ if $0 \leq e[i](k)$. Then, due to Eq. (\ref{eq:w_to_W}), it is sufficient to show that there exists a $c_1$ such that $0 < c_1 \leq W[i,j](k)$ if $0 \leq e[j](k)$.

    (Proof of Condition \ref{con:pos_w}\ref{item:c1a}) 
    Since $ s_1[i,j] > 0 $ and it is not time-varying, there exists a $c_1'$ such that $0 < c_1' \leq s_1[i,j]$. The inequality $0 \leq e[j](k)$ implies that $s_3[i,j] = 1$ for all $i \in \mathcal{I}_1^{m}$ due to Line 6 of the Algorithm \ref{alg:DSMC}.
    Note that
    \begin{equation*}
        \begin{aligned}
            \sum_{i \in \mathcal{I}_1^{m}} E[i,j](k) 
            &= \sum_{i \in \mathcal{I}_1^{m}} \bigg( 0, \Big( e[j](k) - e[i](k) \Big)  A_{a}[i,j] \bigg) \\
            &\leq \sum_{i \in \mathcal{I}_1^{m}} |e[j](k)-e[i](k)|\\
            &\leq (m - 1) |e[j](k)| + \sum_{ \substack{i \in \mathcal{I}_1^{m} \setminus \{ j \} } } |e[i](k)|.
        \end{aligned}
    \end{equation*}
    Since $||e[i](k)||_{\infty} \leq 1$, $\sum_{i  \in \mathcal{I}_1^{m} } |e[i](k)| \leq 2$, and $m \geq 2$, then
    \begin{equation*}
        \begin{aligned}
            \sum_{i \in \mathcal{I}_1^{m}} E[i,j](k) 
            &\leq (m - 1) |e[j](k)| + (2 - |e[j](k)|)\\
            &= m|e[j](k)| + 2 - 2 |e[j](k)|\\
            &\leq m.
        \end{aligned}
    \end{equation*}
    It then follows that Line 5 of the Algorithm \ref{alg:DSMC}
    implies 
    $\frac{x[j](k)}{m} = \min \Big(\frac{x[j](k)}{m}, 1\Big) \leq \min\Big(\frac{x[j](k)}{\sum_{i \in \mathcal{I}_1^{m}} E[i,j](k)}, 1\Big)  \leq s_2[j](k)$. 
    Let $v_{\text{min}} = \min_{i \in \mathcal{I}_1^{m}} v[i]$. Then, the inequality $0 \leq e[j](k)$ implies that $v_{\text{min}} \leq x[j](k)$. 
    Thus, $v_{\text{min}}/m \leq x[j](k)/m \leq s_2[j](k)$.
    Since $v_{\text{min}}/ m$ is not time-varying, there exists a $c_1''$ such that $0< c_1'' \leq v_{\text{min}} / m \leq s_2[j](k)$. Hence, if $0 \leq e[j](k)$, then there exist a $c_1$ such that $0 < c_1 \leq \min\Big(s_1[i,j] , s_2[j](k), s_3[i,j] \Big) = W[i,j](k)$, where $c_1 = \min (c_1', c_1'', 1)$.
    
    (Proof of Condition \ref{con:pos_w}\ref{item:c1b}) 
    Condition \ref{con:pos_w}\ref{item:c1b} requires that $w_{i,j} = 0$ if $\{ i, j \} \in I_{{n_p},n}$, where $I_{{n_p},n}$ is the set defined in Definition \ref{def:erridxset}.
    Note that $\{ i, j \} \in I_{{n_p},n}$ if and only if $e[\ell] < 0$ for all $\ell \in \{ i, j \}$ and $e[l] \geq 0$ for some $l \in I_{ij}$.
    Therefore, if $\{ i, j \} \in I_{{n_p},n}$, both $s_3[i,j](k)$ and $s_3[j,i](k)$ are set to $0$, which imply that both $W[i,j](k)$ and $W[j,i](k)$ are $0$. Hence, $w_{i,j}(k) = 0$ if $\{ i, j \} \in I_{{n_p},n}$.
\end{proof}

The following theorem shows that Algorithm \ref{alg:DSMC} achieves the desired convergence.
\begin{theorem} \label{thm:dsmc_conv}
The probability distribution $x(k)$ exponentially converges to the desired distribution $v$ if the corresponding Markov matrix is synthesized by Algorithm \ref{alg:DSMC}.
\end{theorem}

\begin{proof}
    In Proposition \ref{pro:error_evol}, it is proven that the dynamics of the error vector in Algorithm \ref{alg:DSMC} are identical to the dynamics of the value vector in Theorem \ref{thm:graph_general_2}. Condition \ref{con:pos_w} is used in Theorem \ref{thm:graph_general_2} to prove that the value vector exponentially converges to $\bm{0}$. In Proposition \ref{pro:ams_satis}, it is proven that Algorithm \ref{alg:DSMC} satisfies Condition \ref{con:pos_w} which means that probability distribution $x(k)$ exponentially converges to the desired distribution $v$.
\end{proof}

\vspace{-0.20cm}
\subsection{The Modified DSMC Algorithm} \label{sec:dsmc_2}
In this section, we introduce a shortest-path algorithm that is proposed as a modification to the Metropolis-Hastings algorithm in \cite[Section V-E]{accikmecse2012markov} and integrated with the Markov chain synthesis methods described in \cite{bandyopadhyay2017probabilistic} and \cite{jang2018local}. This algorithm can also be integrated with the DSMC algorithm to further increase the convergence rate. Our approach begins by categorizing the states of the desired distribution.

\vspace{-0.20cm}
\begin{definition} \label{def:recur_trans}
    \textit{(Recurrent and transient states)} 
    The states with non-zero elements in the desired distribution $v$ are called recurrent states. The other states with zero elements in the desired distribution $v$ are called transient states.
\end{definition}

\vspace{-0.20cm}
The shortest-path algorithm is used for the synthesis of the Markov chain for the transient states while the Markov chain for the recurrent states is synthesized by the DSMC algorithm.
Let us split the desired distribution and the Markov matrix to synthesize the Markov matrix for recurrent and transient states separately. Assume that there are $m_r$ recurrent states and $m_t = m - m_r$ transient states of the desired distribution. Then the Markov matrix and the desired distribution are split as
\begin{equation} \label{eq:split}
    \begin{aligned}
        v &= (v_t, v_r), \quad
        M(k) =
        \begin{bmatrix}
            M_1 & \bm{0}\\
            M_2 & M_3(k)\\
        \end{bmatrix}_{m \times m},
    \end{aligned}
\end{equation}
where $v_t = \bm{0}$, $v_t \in \mathbb{R}^{m_t}$, $v_r > \bm{0}$, $v_r \in \mathbb{R}^{m_r}$, $M_1 \in \mathbb{R}^{m_t \times m_t}$ and $M_3(k) \in \mathbb{R}^{m_r \times m_r}$. 
Since the desired distribution and the Markov matrix are partitioned by renumbering the elements, the probability distribution and the adjacency matrix are also partitioned using the same renumbering as
\begin{equation*}
    \begin{aligned}
        x(k) &= (x_{t}(k), x_{r}(k)), \quad
        A_a =
        \begin{bmatrix}
            A_{a_1} & \bm{0}\\
            A_{a_2} & A_{a_3}\\
        \end{bmatrix}_{m \times m}. \\
    \end{aligned}
\end{equation*}

Note that while the time-invariant matrices $M_1$ and $M_2$ in Eq. (\ref{eq:split}) are synthesized by the shortest-path algorithm, the time-varying Markov matrix $M_3(k)$ is synthesized by the DSMC algorithm.
If the shortest-path is to be combined with another Markov chain synthesis algorithm,  it is also necessary to assume that the recurrent states are connected among themselves, which requires the following assumption.
\begin{assumption} \label{asm:strong}
    The graph defined by the recurrent states is connected,
    which means there exists a path between all recurrent states, or equivalently, $(I + A_{a_3})^{m_r-1} > \bm{0}$ for $A_{a_3} \in \mathbb{R}^{m_r \times m_r}$ \cite[section 2.1]{mesbahi2010graph}.
\end{assumption}
Furthermore, we demonstrate that the condition $\bm{1}^T x_{r}(k) = 1$, which is crucial for the Markov chain synthesis algorithm of the recurrent states, is satisfied for all $k \geq K$, $K \in \mathbb{Z}_+$ in the shortest-path algorithm.

Let us first define the following index sets to present the algorithm.
\begin{definition} \label{def:recurridxset}
\textit{(Index sets with respect to recurrent states)}
The index sets $I_r$ and $I_t$ contain the indices of the recurrent and transient states, respectively. Indices of transient states can be split into subsets as $I_s$, $I_{s-1},...$, and so on. The index set $I_{s}$ consists of the indices of the states that are adjacent to the states corresponding to $I_r$ and $I_{s} \cap I_r = \emptyset $, the index set $I_{s-1}$ consists of the indices of states that are adjacent to the states corresponding to $I_{s}$ and $I_{s-1} \cap (I_r \cup I_{s}) = \emptyset$, and so on.
\end{definition}
 
The synthesis of the Markov matrix using the DSMC algorithm and the shortest-path algorithm for the recurrent and transient states of the desired distribution is given as 
\begin{dmath*} 
  M[i,j] =
  \begin{cases}
    \text{as in Algorithm \ref{alg:DSMC}}                     & \text{if $i \in I_r$, $j \in I_r$} \\
    0                                                               & \text{\parbox{15em}{if $i \in I_{s-n}$, $j \in I_r$ for $n \in \mathbb{Z}_+$ } } \\
    $$1 / (\sum_{l \in I_r    } A_a[j,l]$$)& \text{if $i \in I_r$, $j \in I_s$, $A_a[j,i] = 1$} \\
    $$1 / (\sum_{l \in I_{s-n}} A_a[j,l]$$)& \text{\parbox{15em}{if $i \in I_{s-n}$, $j \in I_{s-n-1}$, $A_a[j,i] = 1$ for $n \in \mathbb{Z}_+$}} \\
  \end{cases}
\end{dmath*}

In the proposed shortest-path algorithm, the transition probability from any state in set $I_{s-k-1}$ to any state in set $I_{s-k}$, $k \in \mathbb{Z}_+$, and from any state in set $I_s$ to any state in set $I_r$ is $1$.
Therefore, the total probability of the transient states becomes zero in a finite time.
In \cite{accikmecse2012markov}, it is shown that the condition $\rho(M_1) < 1$ is satisfied using the properties of M-matrices, which are shown in Theorem 2.5.3 (parts 2.5.3.2 and 2.5.3.12) of \cite{horn1994topics}.
\section{An application of the DSMC Algorithm \\ on Probabilistic Swarm Guidance} \label{sec:background}

In this section, we apply the DSMC algorithm to the probabilistic swarm guidance problem and provide numerical simulations that show the convergence rate of the DSMC algorithm is considerably faster as compared to the previous Markov chain synthesis algorithms in \cite{accikmecse2012markov} and \cite{bandyopadhyay2017probabilistic}.

\vspace{-0.20cm}
\subsection{Probabilistic Swarm Guidance} \label{sec:psg_}
Most of the underlying definitions and baseline algorithms in this section are based on \cite{accikmecse2012markov}.
We briefly review this material for completeness.
\vspace{-0.20cm}
\begin{definition} \label{def:bins}
\textit{(Bins)} The operational region is denoted as $\mathcal{R}$. The region is assumed to be partitioned as the union of $m$ disjoint regions, which are referred to as bins $R_i$, $i = 1,...,m$, such that  $\mathcal{R} = \cup_{i=1}^{m} R_i$, and $R_i \cap R_j = \varnothing$ for $i \neq j$. Each bin can contain multiple agents.
\end{definition}

Bins, as in the states of a Markov chain,  in the operational region are represented as vertices of a graph. The allowed transitions between bins are specified by the adjacency matrix of the graph. 
\vspace{-0.20cm}
\begin{definition} \label{def:motcons}
\textit{(Transition constraints)} An adjacency matrix $A_a \in \mathbb{R}^{m \times m}$ is used to restrict the allowed transitions between bins. $A_a[i,j] = 1$ if the transition from $R_i$ to $R_j$ is allowable, and is $0$ otherwise. 
The bins $R_i$ and $R_j$ are called neighboring bins if $A_a[i, j] = 1$.
\end{definition}

\vspace{-0.20cm}
\begin{assumption} \label{asm:strongly_connected}
    The graph describing the bin connections is connected, that is, 
    $(I + A_{a})^{m-1} > \bm{0}$ for $A_{a} \in \mathbb{R}^{m \times m}$ \cite[section 2.1]{mesbahi2010graph}.
\end{assumption}

\vspace{-0.20cm}
Consider a swarm comprised of $N$ agents.
Let $n[i](k)$ be the number of agents in bin $R_i$ at time $k$. Then $x(k) = n(k) / N$ represents the swarm density distribution at time $k$ where $x(k) \geq \bm{0}$ and $\bm{1}^T x(k) = 1$.
It is desired to guide the swarm density distribution to a {\em desired steady-state distribution}
 $v \in \mathbb{R}^m$,  where  $v \geq \bm{0}$ and $\bm{1}^T v = 1$.
The distance between the swarm density distribution and the desired distribution is monitored via the total variation between the swarm density distribution at time $k$ and the desired distribution
\begin{equation*} \label{eq:tv}
    T_{x, v}(k) = \frac{1}{2} ||x(k) - v||_1.
\end{equation*}

\begin{assumption}
All agents can determine their current bins, and they can use this information to move between bins.
\end{assumption}

Each swarm agent changes its bin at time $k$ by using a column stochastic, Markov, matrix $M(k) \in \mathbb{R}^{m \times m}$ called a Markov matrix \cite{horn2012matrix}.
The entries of matrix $M(k)$  define the transition probabilities between bins, that is, for any $k \in \mathbb{Z}_+$ and $i, j \in \mathcal{I}_1^m$,
\begin{equation*}
    M[i,j](k) = \mathcal{P}(r(k+1) \in R_i | r(k) \in R_j ), \\
\end{equation*} 
\textit{i.e.}, an agent in $R_j$ moves to $R_i$ with probability $M[i, j](k)$ at time $k$. 
Algorithm \ref{alg:pga} is an implementation of this probabilistic guidance algorithm. The first step of the algorithm is to determine the agent’s current bin. In the following steps, each agent samples from a uniform discrete distribution and goes to another bin depending on the column of the Markov matrix, which is determined by the agent’s current bin. 
\begin{algorithm}
    \begin{algorithmic}[1]
        \State Each agent determines its current bin, e.g., $r_l(k) \in R_i$.
        \State Each agent generates a random number $z$ that is uniformly distributed in $[0, 1]$.
        \State Each agent goes to bin $j$, \textit{i.e.}, $r_l(k + 1) \in R_j$, if
        $
      \begin{cases}
          $$ \sum_{s=1}^{j-1} M[s,i](k) \leq z \leq \sum_{s=1}^{j} M[s,i](k)$$.   \\
      \end{cases}
        $
    \end{algorithmic}
\caption{Probabilistic Guidance Algorithm \cite{accikmecse2012markov}}
\label{alg:pga}
\end{algorithm}

The main idea of the probabilistic swarm guidance is to drive the propagation of density distribution vector $x(k)$, instead of individual agent positions $\{r_l(k)\}_{l=1}^N$. 
Swarms are viewed as a statistical ensemble of agents to facilitate the swarm guidance problem. The density distribution of the swarm $x(k)$ satisfies the properties of a probability distribution, which are $x(k) \geq \bm{0}$, and $\bm{1}^T x(k) = 1$. However,  the number of agents is finite
and, using the law of large numbers \cite[Section V]{bertsekas2008introduction}, the Algorithm \ref{alg:pga} ensures that $x(k\!~+~\!1)~\!=\!~M(k)~x(k)$ is satisfied
as $N \xrightarrow[]{} \infty$. 
Furthermore, 
any given desired distribution $v \in \mathbb{R}^m$ cannot always be accurately represented by the finite number of agents
of the swarm. 
For example, if $v = [0.5, 0.5]^T$ and $N = 1$, then the best-quantized representation of the distribution that can be obtained is either $[1, 0]^T$ or $[0, 1]^T$. In this case, the total variation between the density distribution of the swarm 
and the desired distribution cannot be less than $0.5$. 
\vspace{-0.20cm}
\begin{definition} \label{def:qe}
    \textit{(Quantization error)} 
    The quantization error is the minimum total variation value that can be achieved by the density distribution of the swarm $x(k) = n(k)/N$
    to represent a given desired distribution $v$. It is denoted as $q_N(v)$ and is the solution to the following problem
    \begin{equation*}
        q_N(v) = \min_n \bigg \{ \frac{1}{2} || n/N - v ||_1 : n \in \mathbb{Z}_+^m, \bm{1}^T n = N \bigg \}.
    \end{equation*}
\end{definition}

Let $Q_N$ be the maximum value of the quantization error that considers the worst possible desired distribution $v$ for a given number of agents $N$ such that 
\begin{equation*}
    Q_N = \max_v \{ q_N(v) : v \in \mathbb{R}_+^m, \bm{1}^T v = 1 \}.
\end{equation*}

The following theorem provides the maximum quantization error value for a given number of agents $N$.
\vspace{-0.20cm}
\begin{lemma} \label{thm:quant}
    For any given desired distribution $v \in \mathbb{R}_+^m$, $\bm{1}^T v = 1 $, the maximum quantization error value between $x(k) = n(k)/N$ and $v$ is not greater than $m/(4N)$, where  $N$ is the total number of agents.
\end{lemma}

\vspace{-0.20cm}
\begin{proof}
    Let us split the integer and fractional part of $Nv$ such that $Nv = t + f$, where $t \in \mathbb{Z}_{+}^m$ and $ \bm{0} \leq f < \bm{1}$. Let $s = \bm{1}^T f$. Note that $s \in \mathbb{Z}_+ $ and $s = N - \bm{1}^T t$. Let $n^*[i]$ be the optimal value of $n[i]$ that minimizes $ |n[i] - Nv[i]| $. Then, $n^*[i] = t[i] $ or $n^*[i] = t[i] + 1$.
    Without loss of generality, suppose that $f[1] \geq f[2] \geq \dots \geq f[m]$. Then      $n^*[i] = 
        \begin{cases}
            t[i] + 1 & \text{if } i \leq s \\
            t[i]     & \text{otherwise } \\
        \end{cases},
        $
    where $\bm{1}^T n = \bm{1}^T t + s = N$.
    Therefore $ || n - vN ||_1 = (1 - f[1]) + \dots + (1 - f[k]) + f[k+1] + \dots + f[m]$. Since $s = \bm{1}^T f$,  $ || n - vN ||_1 = 2 ( f[k+1] + \dots + f[m] )$. Also, $f[k+1] \leq s/m$ since $f[i] \geq f[i+1]$ for all $i \in \mathcal{I}_1^m$ and $\bm{1}^T f = s$. Therefore, $ || n - vN ||_1 \leq 2 (m - s) (s/m) $. The maximum value of $2 (m - s) (s/m) $ is $m/4$, where $s= m/2$. Hence, $|| n - vN ||_1 \leq m/2$ and $\frac{1}{2} || n/N - v ||_1 \leq m/(4N)$.
\end{proof}

\subsection {Synthesis of the Markov Matrix} 
The DSMC algorithm presented in Section \ref{sec:dsmc_1} is employed for synthesizing the Markov chain, which serves as the stochastic policy for the swarm agents.
It is worth noting that the bins comprising the operational region, as defined in Definition \ref{def:bins}, determine the vertices of the uniform graph in Definition \ref{def:graph_def}. Consequently, these vertices correspond to the states of the Markov chain defined in Definition \ref{def:Markov}. Similarly, the transition constraints of the swarm, defined by an adjacency matrix in Definition \ref{def:motcons}, determine the adjacency matrix of the uniform graph in Definition \ref{def:graph_def}, which subsequently corresponds to the adjacency matrix of the Markov chain in Definition \ref{def:Markov}.
The connectivity requirement for the vertices of the consensus protocol, as outlined in Assumption \ref{asm:str_con_graph}, is satisfied by Assumption \ref{asm:strongly_connected}. Furthermore, if the recurrent states of the desired distribution are also connected among themselves, then Assumption \ref{asm:strong} is satisfied, allowing the utilization of the Modified DSMC algorithm presented in Section \ref{sec:dsmc_2} for Markov chain synthesis.
Since the DSMC algorithm is state-dependent, we use the density distribution of the swarm for feedback, which requires
the following assumption.
\vspace{-0.20cm}
\begin{assumption} \label{asm:asm_dens}
All agents know the density values of their own and neighboring bins.
\end{assumption}

\vspace{-0.20cm}
A complex communication architecture is not required since communication only with neighboring bins is sufficient for an agent to determine its transition probabilities. 
If agents have only access to the number of agents of their own and neighboring bins, then they also need to know the total number of agents in the swarm, which is global information, to determine their density values. Distributed consensus algorithms, which work under the strongly-connected communication network topology assumption, can be used to estimate the total number of agents of the swarm. The consensus protocol that is used to estimate the density distribution of the swarm in \cite[Remark 13]{bandyopadhyay2017probabilistic} can also be used to estimate the total number of agents of the swarm.
\vspace{-0.3cm}
\subsection{Numerical Simulations} \label{sec:results}
\begin{figure*}[!hbt]
    \centering
    \subfloat[$t=0$\label{fig:mh_imc_su_0}]{%
        \includegraphics[width=0.20\linewidth]{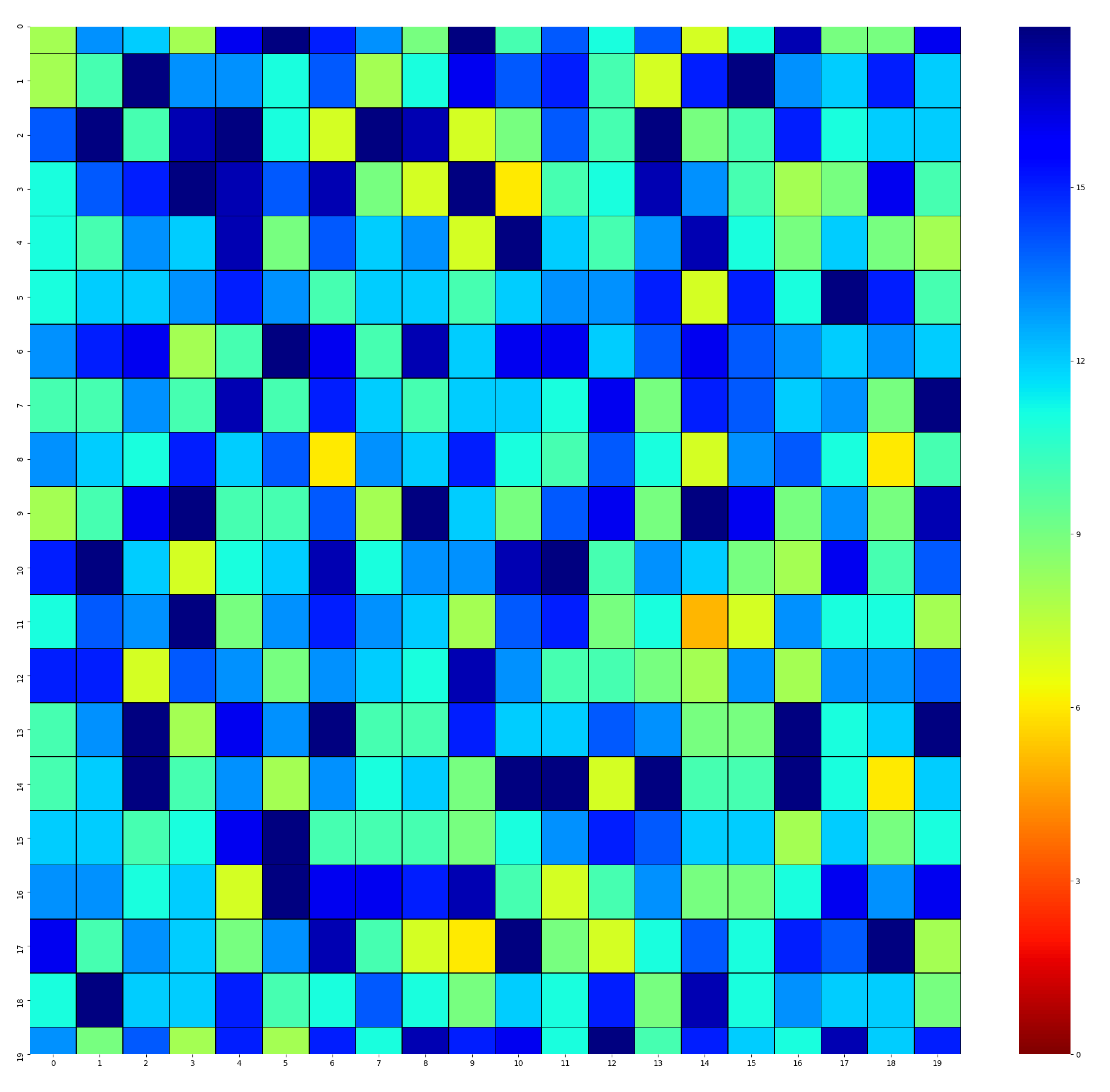}}
    \hfill
    \subfloat[$t=250$\label{fig:mh_imc_su_250}]{%
        \includegraphics[width=0.20\linewidth]{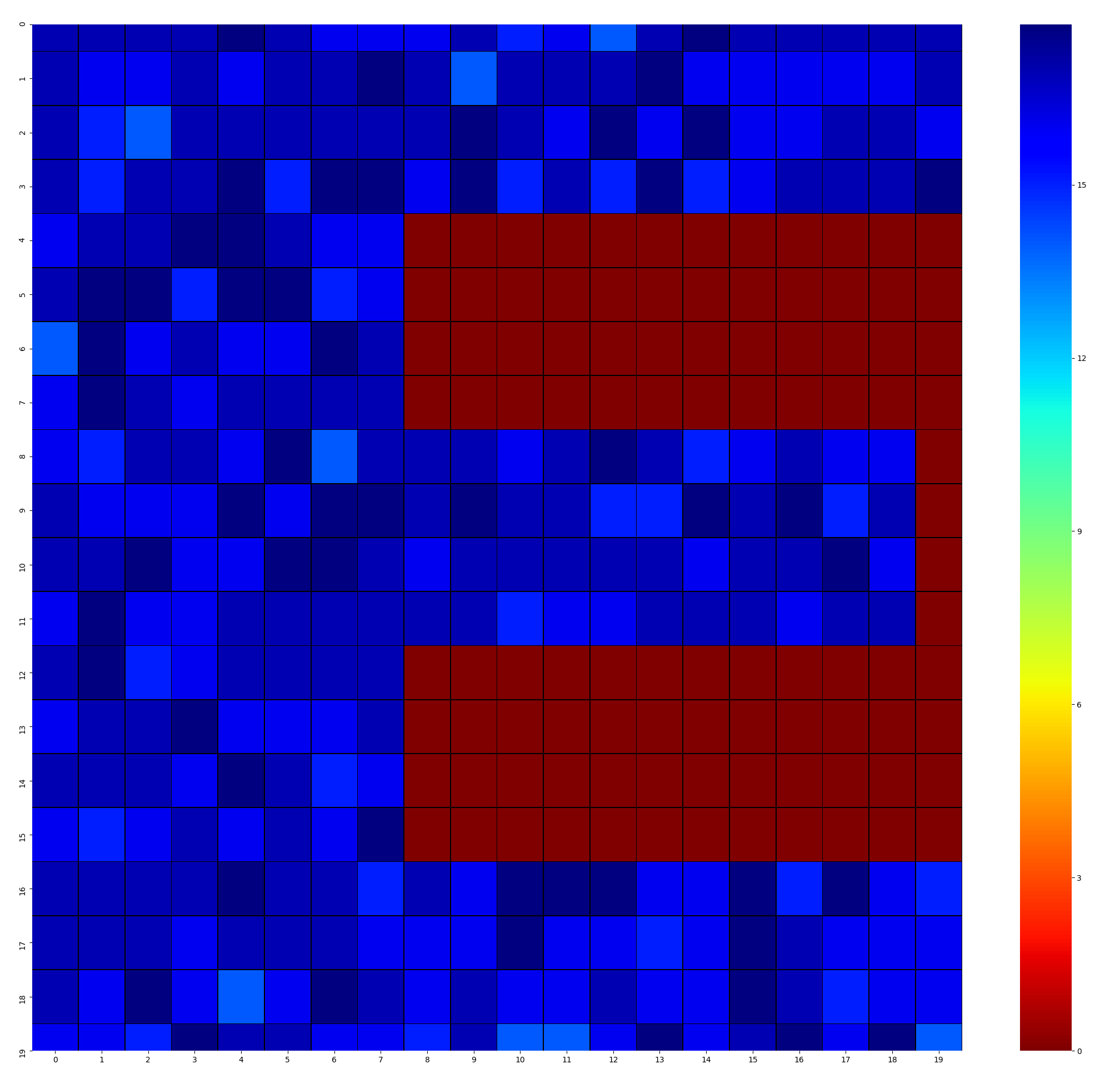}}
    \hfill
    \subfloat[$t=251$\label{fig:mh_imc_su_251}]{%
        \includegraphics[width=0.20\linewidth]{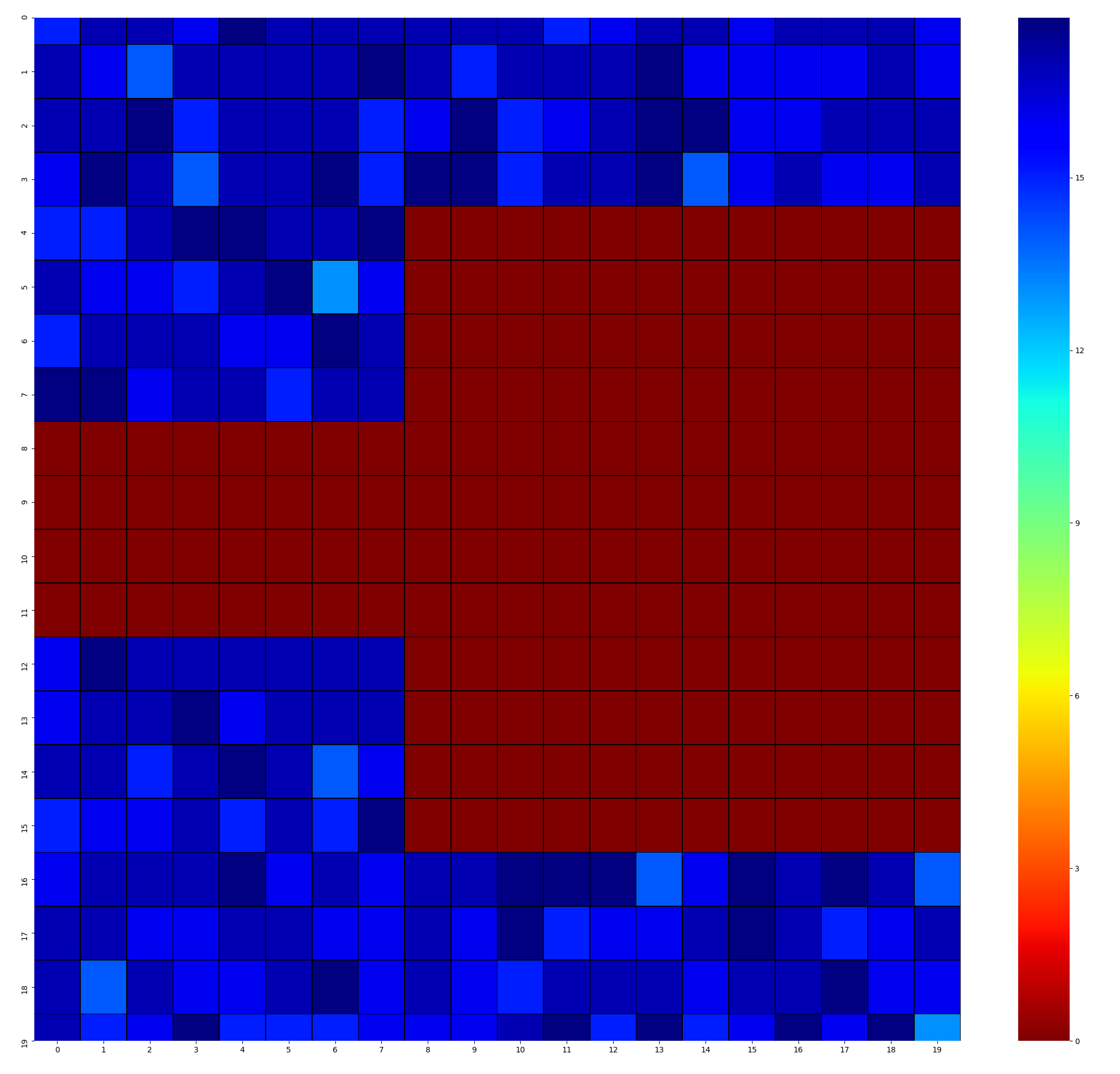}}
    \hfill
    \subfloat[$t=750$\label{fig:mh_imc_su_750}]{%
        \includegraphics[width=0.20\linewidth]{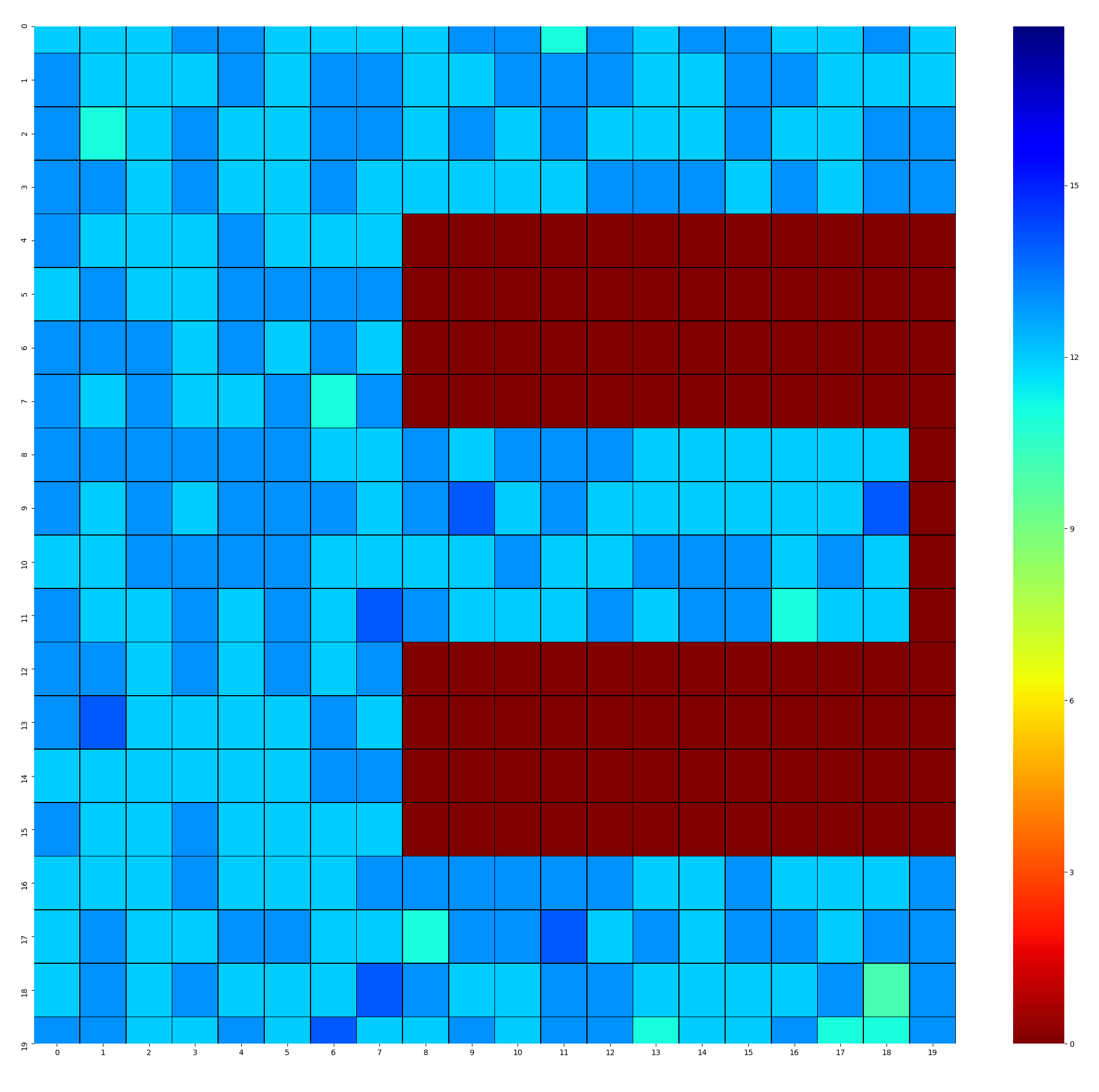}}
    \caption{Representation of the distribution of the swarm for the time-steps $0$, $250$, $251$, and $750$, respectively. There are $400$ ($20 \times 20$) bins and $5000$ agents in the operational region at the beginning of the simulation. The agents converge to the ``E'' letter in $250$ time-steps and approximately $1/3$ agents are removed from the operational space. Then, the remaining agents converge to the desired distribution again in $500$ time-steps.}
    \label{fig:MHNEW}
\end{figure*}

\begin{figure*}[!hbt]
    \centering
    \subfloat[$t=0$\label{fig:g0}]{%
        \includegraphics[width=0.232\linewidth]{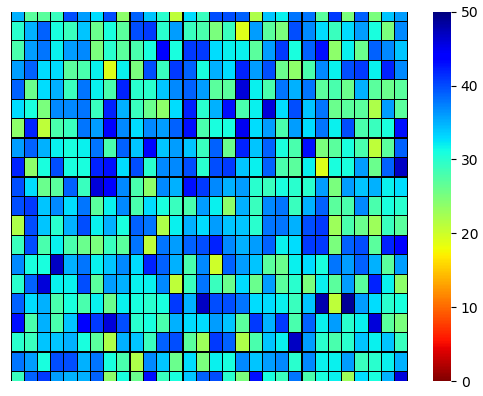}}
    \hfill
    \subfloat[$t=40$\label{fig:g40}]{%
        \includegraphics[width=0.232\linewidth]{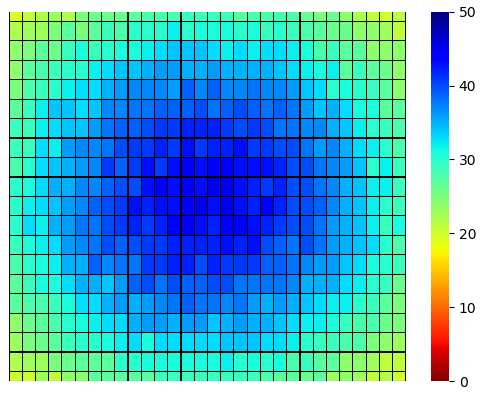}}
    \hfill
    \subfloat[$t=80$\label{g80}]{%
        \includegraphics[width=0.232\linewidth]{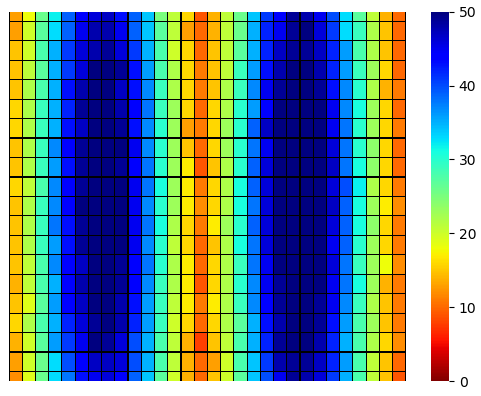}}
    \hfill
    \subfloat[$t=120$\label{fig:g120}]{%
        \includegraphics[width=0.232\linewidth]{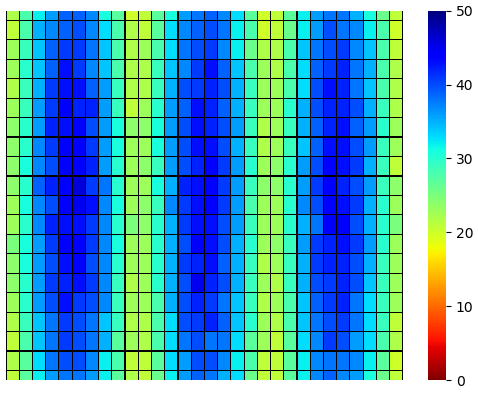}}\\
    \subfloat[$t=160$\label{fig:g160}]{%
        \includegraphics[width=0.232\linewidth]{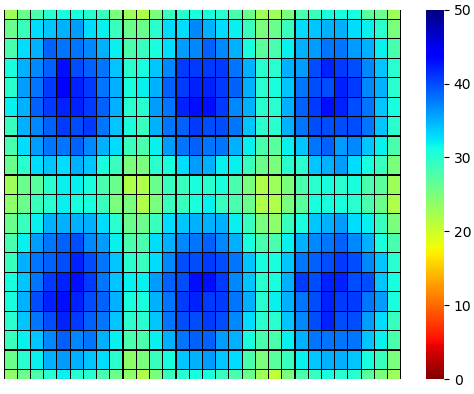}}
    \hfill
    \subfloat[$t=161$\label{fig:g161}]{%
        \includegraphics[width=0.232\linewidth]{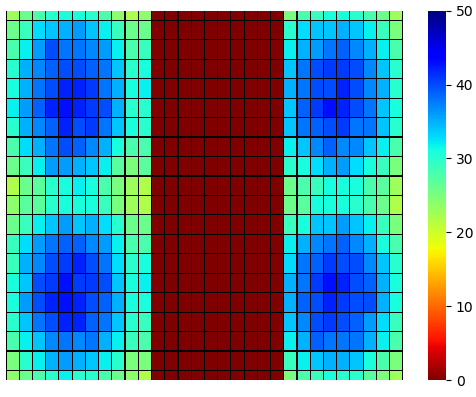}}
    \hfill
    \subfloat[$t=200$\label{fig:g200}]{%
        \includegraphics[width=0.232\linewidth]{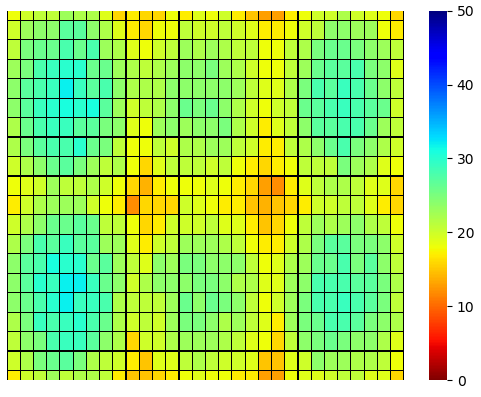}}
    \hfill
    \subfloat[$t=201$\label{fig:g201}]{%
        \includegraphics[width=0.232\linewidth]{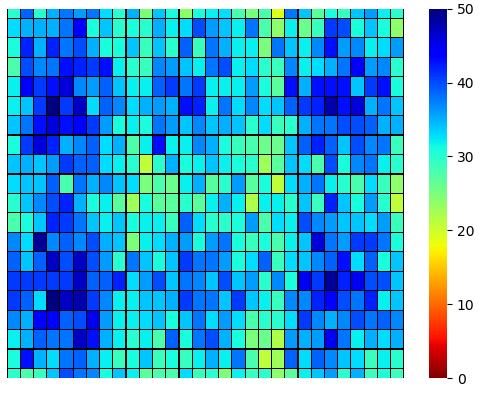}}\\
    \subfloat[$t=240$\label{fig:g240}]{%
        \includegraphics[width=0.232\linewidth]{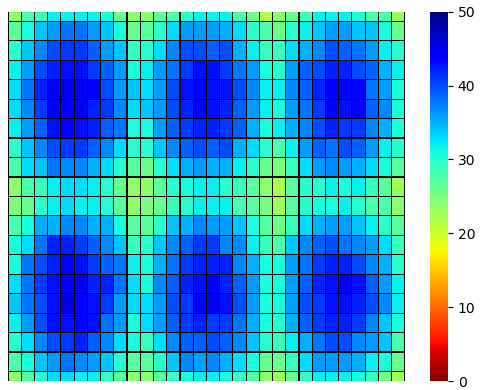}}
    \hfill
    \subfloat[$t=280$\label{fig:g280}]{%
        \includegraphics[width=0.232\linewidth]{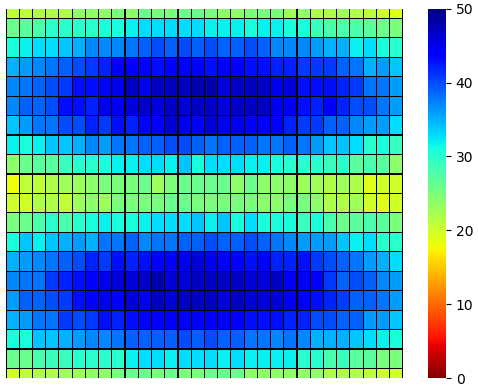}}
    \hfill
    \subfloat[$t=320$\label{fig:g320}]{%
        \includegraphics[width=0.232\linewidth]{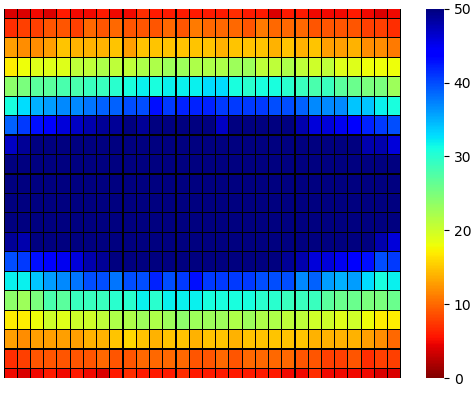}}
    \hfill
    \subfloat[$t=360$\label{fig:g360}]{%
        \includegraphics[width=0.232\linewidth]{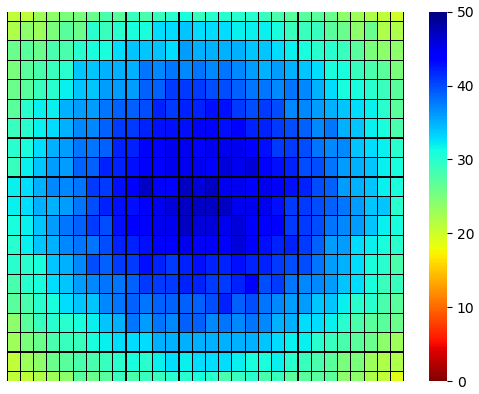}}
    \caption{ Representation of the distribution of the swarm for several time-steps. There are $600$ ($20 \times 30$) bins and $20000$ agents in the operational region at the beginning of the simulation. Swarm distribution converges to a different multimodal Gaussian distribution in every $40$ time-step except the time-steps between $160$ and $240$. Agents converge to the $4$ different multimodal Gaussian distributions up to the $160^{th}$ time-step. At the $161^{th}$ time-step, approximately $1/3$ agents are removed and the remaining agents converge to the same desired distribution in $40$ time-steps. Then, $7500$ agents are uniformly added to the operational region at $201^{th}$ time-step and all agents converge to the same desired distribution again in $40$ time-steps. After the $240^{th}$ time-step, agents converge to the $3$ new multimodal Gaussian distributions up to the $360^{th}$ time-step.}
    \label{fig:gauss}
\end{figure*}

The convergence performance of the DSMC algorithm is demonstrated on the probabilistic swarm guidance application with numerical simulations shown in Figure \ref{fig:MHNEW} and \ref{fig:gauss}.
We monitor the total variation between the density distribution of the swarm and the desired distribution. 
The finite number of agents causes the quantization error outlined in Definition \ref{def:qe}.
Theorem \ref{thm:quant} limits the quantization error to $m/(4N)$, where $m$ is the number of bins and $N$ is the number of agents. 

\begin{figure*}[!hbt]
    \centering
       \includegraphics[width=0.43\linewidth]{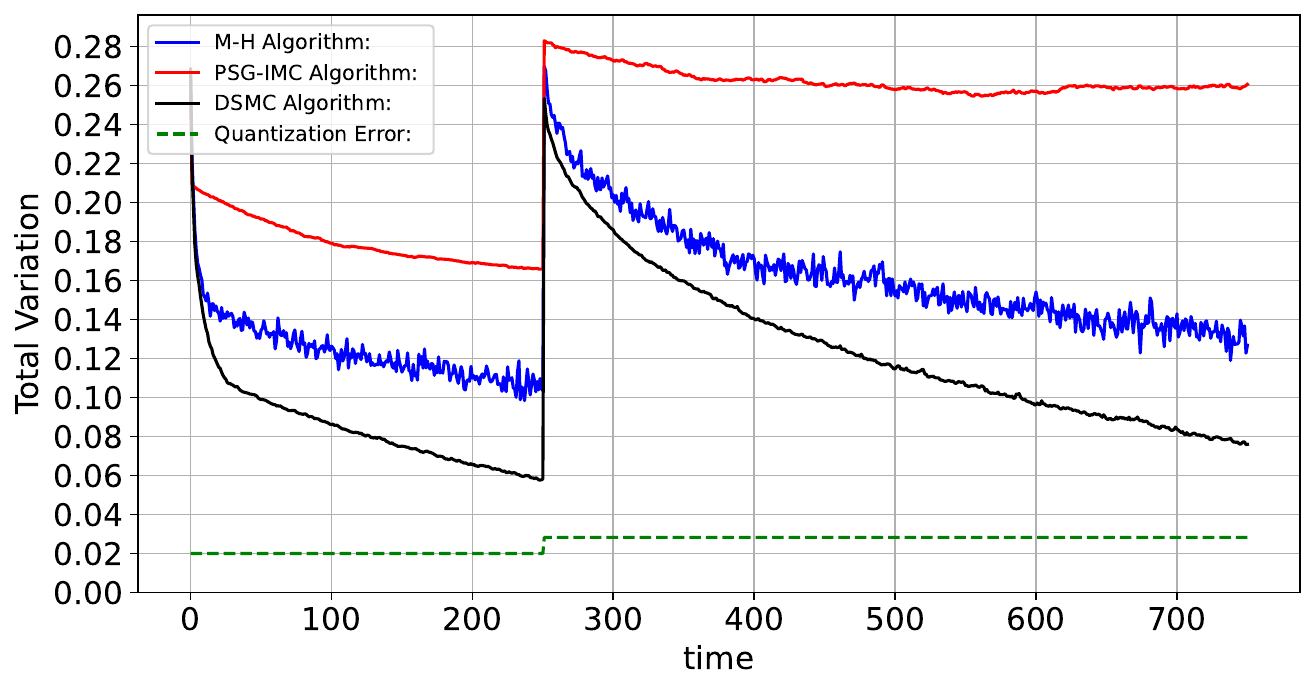}
    \hfill
       \includegraphics[width=0.43\linewidth]{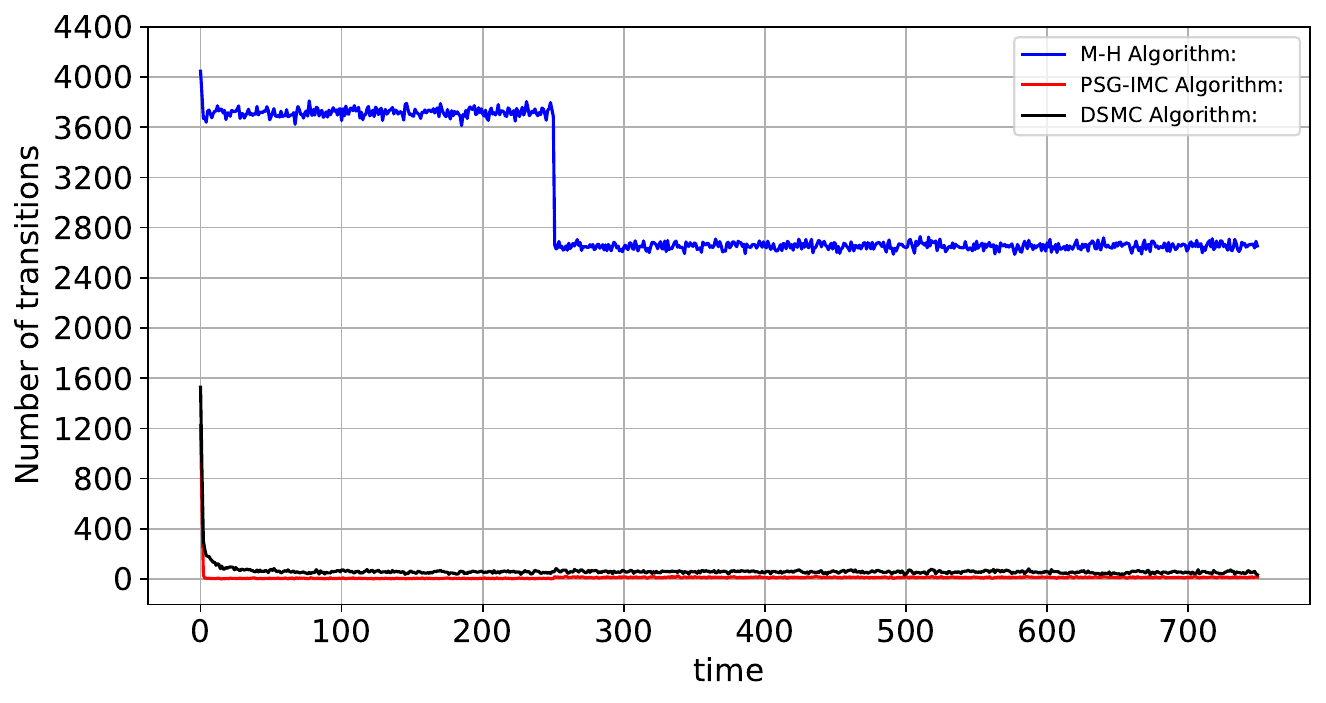}
    \caption{Comparison of change of the total variation and the number of transitions with time for the algorithms. In this case, the adjacency matrix only allows agents to transition to $1$-step away bins.
    }
    \label{fig:MHNEWtotvarntrans}
\end{figure*}

\begin{figure*}[!hbt]
    \centering
        \includegraphics[width=0.43\textwidth]{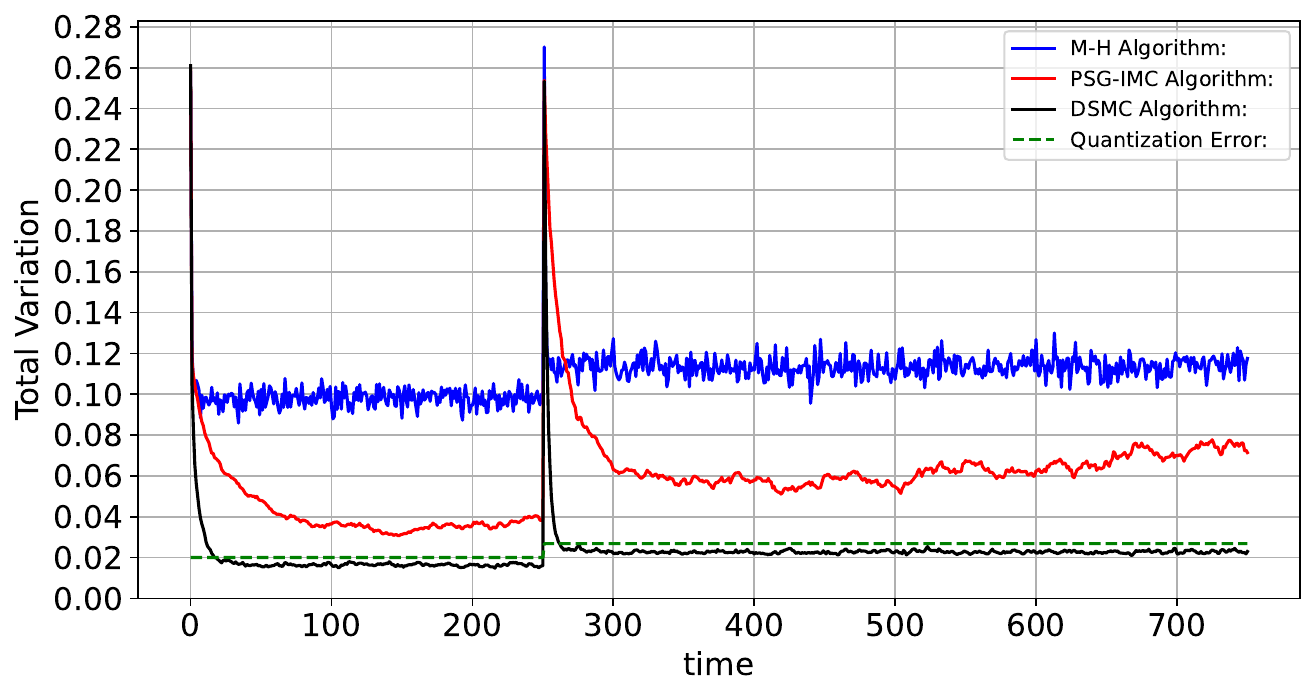}
        \label{fig:10_new_idea_Total_Variation}
    \hfill
        \includegraphics[width=0.43\textwidth]{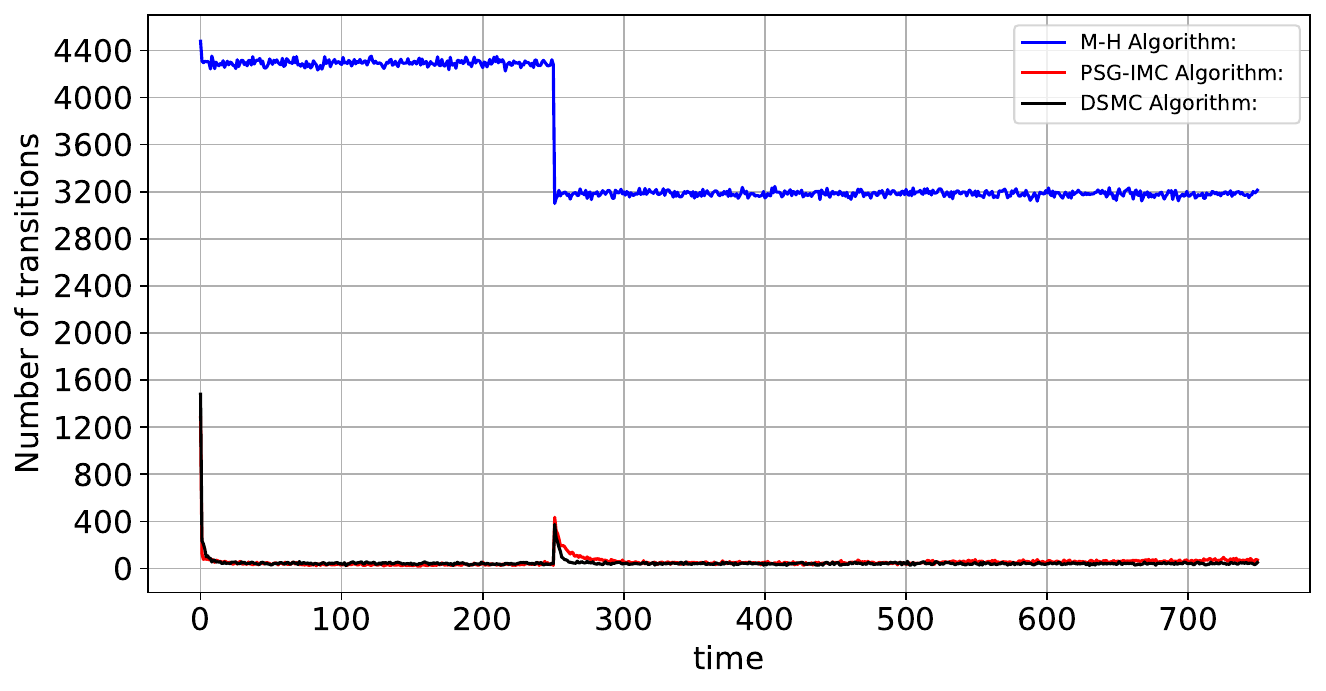}
        \label{fig:10_new_idea_N_of_transitions}
    \caption{Comparison of change of the total variation and the number of transitions with time for the algorithms. In this case, the adjacency matrix allows agents to transition to $10$-step away bins.}
    \label{fig:10_MHNEWtotvarntrans}
\end{figure*}

In the first simulation, which is the same numerical example presented in \cite{accikmecse2012markov}, 
agents converge to the letter ``E" in $250$ time-steps. Then approximately $1/3$ agents are removed and the remaining agents converge to the desired distribution again in $500$ time-steps, demonstrating the ``self-repair" property of the proposed algorithm for swarm guidance. Comparisons of the total variation and the total number of transitions for the M-H, PSG-IMC, and the DSMC
algorithms are given in Figures \ref{fig:MHNEWtotvarntrans}, \ref{fig:10_MHNEWtotvarntrans}, and Tables \ref{table:MHNEWtotvarntrans}, \ref{table:10_MHNEWtotvarntrans} for two different cases. 
In the first case, the adjacency matrix only allows transitions to the bins, which are above, below, right, and left to a given bin, i.e., the bins that are 1-step away. In the second case, the adjacency matrix allows transitions to 10-step away bins. 
When compared to the M-H and PSG-IMC algorithms, in total variation at the $750^{th}$ time-step, the DSMC algorithm provides approximately $1.26$ and $60.67$ times improvement in the speed of convergence (the ratios of the total amount of changing of the total variations in the given time-step) in the first case and $1.51$, and $1.25$ times improvement in the second case. 
The PSG-IMC algorithm does not perform as well
in the first case because of the issues caused by having a sparse adjacency matrix as discussed in Section \ref{sec:related_work}. In the second case, which has a much denser adjacency matrix,
the PSG-IMC algorithm provides approximately $1.21$ times improvement in total variation compared to the M-H algorithm but the DSMC algorithm still provides faster convergence than the PSG-IMC algorithm. 
Furthermore, the total variation value of the DSMC algorithm converges to a value below the maximum value of the quantization error provided by Theorem \ref{thm:quant}.
In the DSMC algorithm, as the distribution converges, the Markov matrix turns into an identity matrix, and unnecessary transitions are avoided as in the PSG-IMC algorithm.
Since the transition of agents rarely occurred in the PSG-IMC algorithm in the first case due to the sparse adjacency matrix, the total number of transitions of the DSMC algorithm is approximately $49.56$ times less than the M-H algorithm, and $4.93$ times more than the PSG-IMC algorithm. For the second case, the total number of transitions of the DSMC algorithm is approximately $74.96$ times less than the M-H algorithm, and $1.17$ times less than the PSG-IMC algorithm. 

\begin{figure*}[!hbt]
    \centering
        \includegraphics[width=0.43\textwidth]{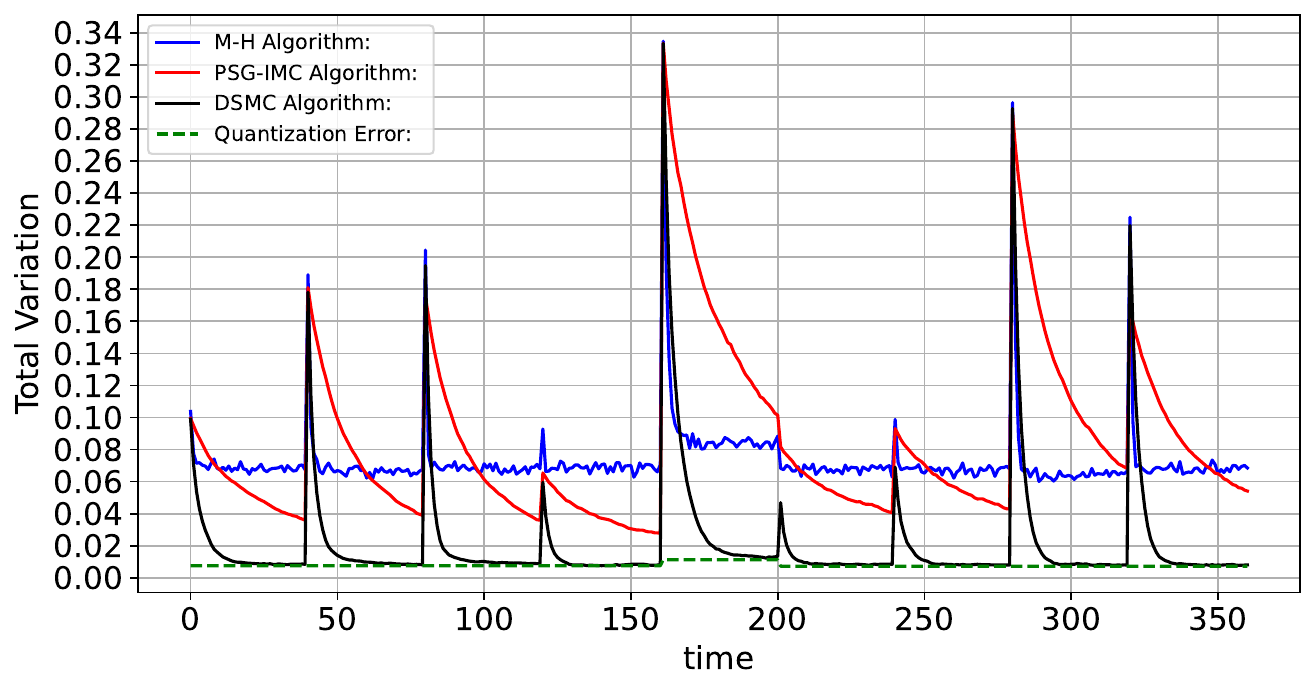}
    \hfill
        \includegraphics[width=0.43\textwidth]{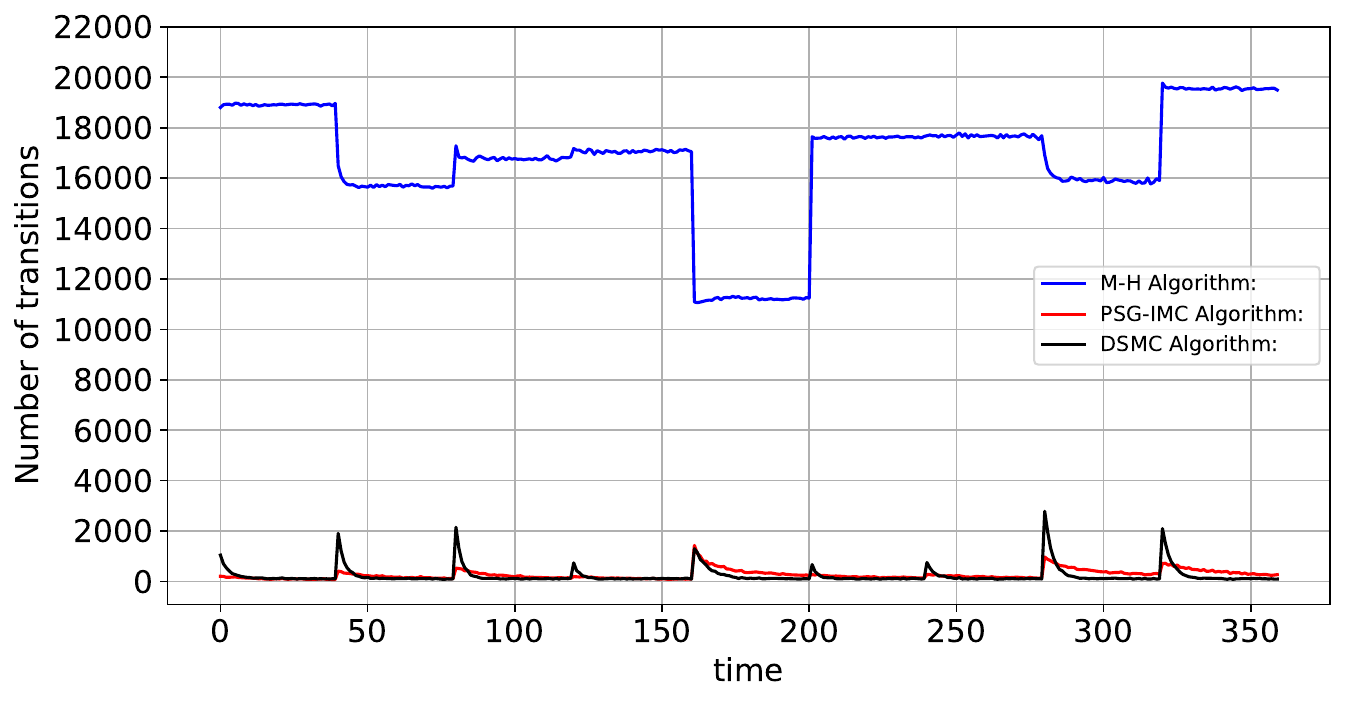}
    \caption{Comparison of change of the total variation and the number of transitions with time for the algorithms. }
    \label{fig:gauss_tvnt}
\end{figure*}

In the second simulation, a more comprehensive comparison is provided. The swarm distribution converges to various multimodal Gaussian distributions (i.e., a mixture of multiple Gaussian distributions). The desired distribution is changed to another multimodal Gaussian distribution every $40$ time-steps, except for the time-steps between $160$ and $240$. Up to the $160^{th}$ time-step, agents converge to the $4$ different multimodal Gaussian distributions. To show the self-repair property of the swarm, approximately $1/3$ of the agents are removed at the $161^{th}$ time-step and the remaining agents converge to the same desired distribution in $40$ time-steps. Moreover, $7500$ agents are uniformly added to the operational region at $201^{th}$ time-step and all agents converge to the same desired distribution again in $40$ time-steps. After the $240^{th}$ time-step, agents converge to the $3$ new multimodal Gaussian distributions up to the $360^{th}$ time-step. Comparisons of the total variation and the total number of transitions for the M-H, PSG-IMC, and DSMC algorithms are given in Figure \ref{fig:gauss_tvnt} and Table \ref{table:gauss_tvnt}. In this case, the adjacency matrix allows transitions to 10-step away bins. 
When compared to the M-H and PSG-IMC algorithms, the DSMC algorithm provides a significant improvement in the speed of convergence for each desired multimodal Gaussian distribution, as well as causing fewer transitions compared to both algorithms. Similar to the first simulation, the total variation value of the DSMC algorithm reaches the maximum value of the quantization error given in Theorem \ref{thm:quant}. Furthermore, when some agents are removed or new agents are added to the operational region, the DSMC algorithm recovers the desired distribution faster than the previous algorithms. 

\begin{table*}[!hbt] 
\caption{Comparison of change of the total variation and the total number of transitions with time for the algorithms. In this case, the adjacency matrix only allows agents to transition to $1$-step away bins.}
\centering
\begin{tabular}{|c|c|c|c|c|c|c|c|c|c|}
\hline
                  & \multicolumn{4}{c|}{Total Variations} & \multicolumn{3}{c|}{Ratios of change of Total Variation}    & \multicolumn{2}{c|}{Total Number of Transitions} \\ \hline
Time              & 0       & 250     & 251     & 750     & 0-250   & 250-750 & 0-750          & 750                   & 750 \textit{(Ratios)}   \\ \hline
M-H Algorithm     & 0.268   & 0.107   & 0.273   & 0.123   & 1.27    & 1.16    & \bm{$1.26$}    & 2282394               & \bm{$49.56$}                   \\ \hline
PSG-IMC Algorithm & 0.260   & 0.169   & 0.282   & 0.257   & 2.25    & 6.96    & \bm{$60.67$}   & 9348                  & \bm{$1/(4.93)$}             \\ \hline
DCMS Algorithm    & 0.261   & 0.056   & 0.253   & 0.079   & 1.00    & 1.00    & \bm{$1.00$}    & 46051                 & \bm{$1.00$}                    \\ \hline
\end{tabular}
\label{table:MHNEWtotvarntrans}
\end{table*}

\begin{table*}[!hbt] 
\caption{Comparison of change of the total variation and the total number of transitions with time for the algorithms. In this case, the adjacency matrix allows agents to transition to $10$-step away bins.}
\centering
\begin{tabular}{|c|c|c|c|c|c|c|c|c|c|}
\hline
                  & \multicolumn{4}{c|}{Total Variations} & \multicolumn{3}{c|}{Ratios of change of Total Variation}  & \multicolumn{2}{c|}{Total Number of Transitions} \\ \hline
Time              & 0       & 250     & 251     & 750     & 0-250   & 250-750 & 0-750           & 750           & 750 \textit{(Ratios)} \\ \hline
M-H Algorithm     & 0.271   & 0.101   & 0.269   & 0.113   & 1.44    & 1.48    & \bm{$1.51$}     & 2663105       & \bm{$74.96$}                 \\ \hline
PSG-IMC Algorithm & 0.263   & 0.039   & 0.253   & 0.072   & 1.09    & 1.28    & \bm{$1.25$}     & 41483         & \bm{$1.17$}                  \\ \hline
DCMS Algorithm    & 0.261   & 0.017   & 0.253   & 0.022   & 1.00    & 1.00    & \bm{$1.00$}     & 35529         & \bm{$1.00$}                  \\ \hline
\end{tabular}
\label{table:10_MHNEWtotvarntrans}
\end{table*}

\begin{table*}[!hbt] 
\caption{ Comparison of change of the total variation and the total number of transitions with time for the algorithms. }
\centering
\begin{tabular}{|c|c|c|c|c|c|c|c|c|c|c|c|c|c|c|c|}
\hline
                  & \multicolumn{9}{c|}{Ratios of change of Total Variation} & \multicolumn{2}{c|}{Total Number of Transitions} \\ \hline
Time              & 0-40   & 40-80  & 80-120  & 120-160   & 161-200 & 201-240   & 240-280   & 280-320   & 320-360   & 360       & 360 \textit{(Ratios)} \\ \hline
M-H Algorithm     & 2.57   & 1.47   & 1.34    & 2.00      & 1.29    & 7.06      & 1.77      & 1.24      & 1.34      & 6058502   & 77.09     \\ \hline
PSG-IMC Algorithm & 1.37   & 1.24   & 1.36    & 1.36      & 1.38    & 0.94      & 1.26      & 1.25      & 1.89      & 95277     & 1.21  \\ \hline
DCMS Algorithm    & 1.00   & 1.00   & 1.00    & 1.00      & 1.00    & 1.00      & 1.00      & 1.00      & 1.00      & 78592     & 1.00      \\ \hline
\end{tabular}
\label{table:gauss_tvnt}
\end{table*}
\section{Conclusion and Future Works}\label{sec:conclusion}

This paper introduces a decentralized state-dependent Markov chain synthesis (DSMC) algorithm with an application to the probabilistic swarm guidance problem.
The proposed algorithm is theoretically proven to exhibit exponential convergence, and numerical experiments confirm faster convergence compared to existing homogeneous and time-inhomogeneous Markov chain synthesis algorithms, which respectively guarantee exponential and asymptotic convergence. 
Furthermore, it is observed that the number of state transitions is relatively low for the fast convergence rates it provides when compared with existing algorithms.

For the probabilistic swarm guidance application, removing the assumption that agents have access to density values of their own and neighboring bins will be the subject of future studies. A useful extension of this research may involve imposing safety constraints on the density distribution of the swarm, such as density upper bounds or density rate bounds. Additional future work may include the swarm engagement problem, which considers matching the distribution of a non-collaborative swarm whose density distribution evolves with a known Markov chain. In that case, having a fast converging algorithm, such as DSMC, would possibly be advantageous to quickly react to such changing desired distributions. 
\section{Acknowledgement} \label{sec:acknowledgement}
The authors would like to thank the anonymous reviewers for their helpful and constructive comments that greatly contributed to improving the final version of the paper.

\bibliographystyle{ieeetr}
\bibliography{root}
\begin{IEEEbiography}[{\includegraphics[width=1in,height=1.25in,clip,keepaspectratio]{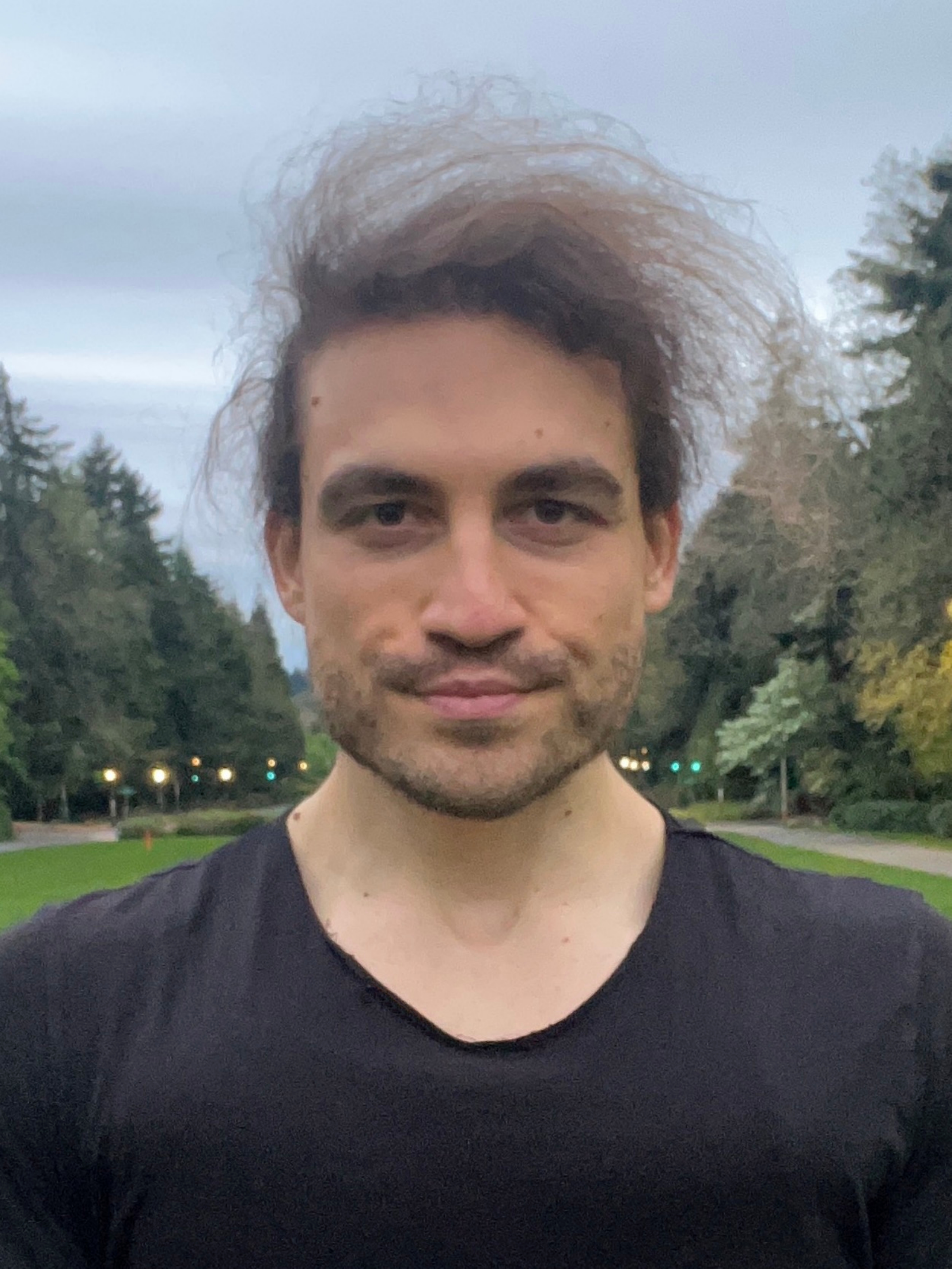}}]{Samet Uzun} 
(Student Member, IEEE) received the B.Sc. and M.Sc. degrees in Department of Aeronautics and Astronautics from Istanbul Technical University, Istanbul, Turkey, in 2018 and 2020, respectively. He is currently pursuing a Ph.D. degree in the Department of Aeronautics and Astronautics from the University of Washington, Seattle, WA, USA. His research interests include optimization, control, and machine learning.
\end{IEEEbiography}

\begin{IEEEbiography}[{\includegraphics[width=1in,height=1.25in,clip,keepaspectratio]{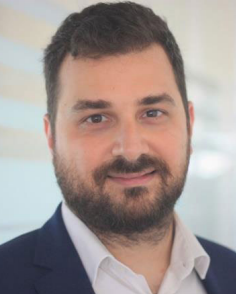}}]{N. Kemal Üre} 
(Member, IEEE) received the B.Sc. and M.Sc. degrees in Department of Aeronautics and Astronautics from Istanbul Technical University, Istanbul, Türkiye, in 2008 and 2010, respectively, and the Ph.D. degree in Department of Aeronautics and Astronautics from the Massachusetts Institute of Technology, Cambridge, MA, USA, in 2015. He is currently an Associate Professor with the Department of Artificial Intelligence and Data Engineering, Istanbul Technical University. His main research interests include applications of deep learning and deep reinforcement learning for autonomous systems, large scale optimization, and the development of high-performance guidance navigation and control algorithms.
\end{IEEEbiography}

\begin{IEEEbiography}[{\includegraphics[width=1in,height=1.25in,clip,keepaspectratio]{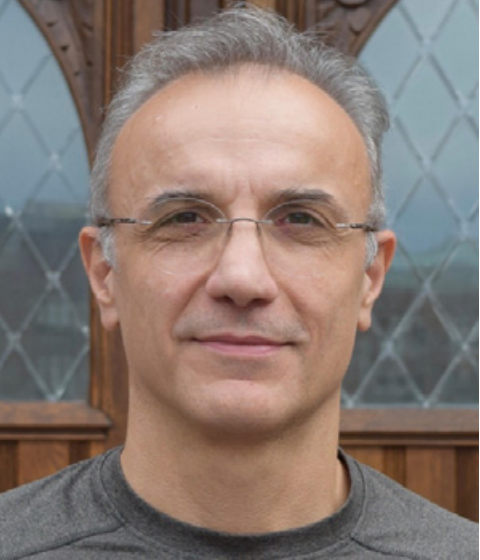}}]{Beh{\c{c}}et A{\c{c}}{\i}kme{\c{s}}e} 
(Fellow, IEEE) received the Ph.D. degree in aerospace engineering from Purdue University, West Lafayette, IN, USA, in 2002.
Previously, he was a Senior Technologist with JPL and a Faculty Member with the University of Texas at Austin, Austin, TX, USA. At JPL, he developed flyaway control algorithms that were successfully used in the landing of Curiosity and Perseverance rovers on Mars. He is currently a Professor with the University of Washington,
Seattle, WA, USA. His research interests include robust and nonlinear control, convex optimization and its applications control theory and its aerospace applications, and Markov decision processes.
Dr. A{\c{c}}{\i}kme{\c{s}}e was the recipient of many NASA and JPL achievement awards for his contributions to spacecraft control in planetary landing, formation flying, and asteroid and comet sample return missions, and also the NSF CAREER Award, IEEE CSS Award for Technical Excel- lence in Aerospace Control and IEEE CSM Outstanding paper award in 2023. He is a Fellow of AIAA and IEEE.
\end{IEEEbiography}

\end{document}